\documentclass[a4paper,10pt]{article}
\usepackage[utf8]{inputenc}
\usepackage{amsmath}
\usepackage{amssymb}
\begin{document}
 
\title{On a weak Lie type for vectors}
\author{Tove Dahn}
\maketitle

 \subsection{Introduction}
 
 The theory of Lie transformation groups, is applied to vector valued distributions, to give
 a scalar type for change of local coordinates, as minimal invariant sets. Conjuation is relative convexity.
 Under the condition that $R(G)^{\bot}$ is $L^{p}-$ compact, with an approximation property, we discuss
 continuations $\mathcal{L}_{ac} \rightarrow \mathcal{L}_{c}$. In particular, we discuss non-regular
 continuations of spiral type. We discuss P-convexity for hypoelliptic operators, as a condition on
 principal operators in Riemann geometry.
 \section{Invariants}

 \subsection{Invariant sets}
 Our starting point is $d U \in \mathcal{G}$ according to Lie, further in the sense of functionals, 
 $Uf = \int f d U(x,y)$, that gives the
 action of movement, with $U \in G$.
 An invariant $\Phi$ (\cite{Lie88}) to a transformation $x_{j} \rightarrow f_{i}$, is given by $\Phi(f_{1}(x),\ldots,f_{n}(x))=\Phi(x_{1},\ldots,x_{n})$.
 Given $X(\Phi)$ a tangent, with $X \in \mathcal{G}_{r}$ , that is $\frac{\delta \Phi}{\delta t}=X(\Phi)$,
 we have that $\Phi$ is characteristic for the transformation. Given r=1, the only invariants to $\mathcal{G}_{1}$,
 are $\Phi$, so that $X(\Phi)=0$.
 
 Assume that $\Psi_{k}$, $k=1,\ldots,n-q$, related to a q-parameter system $X_{j}(\Psi_{k})=0$, $j=1,\ldots,q$, we then have 
 $\{ \Psi_{k}=C_{k} \}$ defines a $(n-q)-$ fold decomposition $\tilde{\Omega}$, the characteristic 
 manifold. Starting from a point P fixed in $E_{n}$, a n-dimensional real space and a complete system, 
 $\tilde{\Omega}$ gives a decomposition of $E_{n}$. 
 
 \newtheorem{prel}{Lemma}[section]
 \begin{prel}
 In the case of Radon-Nikodym, $d U(f)=\rho d I (f)=\frac{d}{d t}(\rho f)d t$, for $\rho \in L^{1}(d I)$, we get $\rho = const$ 
 are characteristic sets to the movement. We write $\Omega(d U)=\{ d U = d I \}$.
 Consider $\mathcal{W}$, a Riemann surface, with $\rho \in \mathcal{H}(\mathcal{W})$ and $\frac{d \rho}{d t}=0$ implies $\rho=const$, 
 we then have for r=1, that invariants are given by $\rho=const$ on $\mathcal{W}$. The condition $\frac{d f}{d t}=d U(f)=0$ 
 gives that U is analytic over first surfaces to f according to $\mathcal{W}$.
 \end{prel}
 
A minimal invariant set as above, defines the type for the movement. 
Assume that $X_{j}(\Psi_{k})=0$ corresponding to $d UW f=0$, with $W f = C_{k} f$, we write $Wf-If=0$, that is $UW - W + W -I$,
where $W-I$ is given by n-q transformations. Given $UW - I=0$ over f, then U is invertible over f, 
but we may have $W \notin G$
 
 More precisely, UW-W+W-I= (U-I)W + W-I, that is R(W) according to the above, is a domain for absolute continuity relatively 
 U and in this case, given
 $f \in \Omega(UW)$ (invariant set), we must have $W=I$ over $\Omega(UW)$. Note however, if $W \in G$, 
 we have that $UW=I$ over
 $\Omega(UW)$, that is $U=W^{-1}$ over $\Omega(UW)$. Given W=I over $\Omega(UW)$, the same holds for U.
 
 Starting from $d U=\rho d V$, where $\rho$ is algebraic, we assume existence of $d W \in G$, 
 such that $d U=d WV$. If we assume $R(V)$ has 
 the approximation property, this property can be continued to $R(U)$. For instance, given W absolute continuous 
 over $R(V)$ and
 $d U=d WV=0$ implies $WV-V=0$, that is W is the restriction of U to $R(V)$. If we further assume $V-I=0$,
 we have invariants to U.  
 
 Assume now that we have $L^{p}-$ compactness for $R(V)^{\bot}$.
 More precisely, assume that $d U=\rho d V$, where $\mid \rho' \mid \leq 1$, for $\rho'=\rho -1$, gives a symmetric 
 neighborhood (in G) of invariant sets. Assume that $\rho'$ analytic, that is we start with an analytic polyhedron. 
 In particular, where $\frac{\delta \rho}{\delta x_{j}}=0$, $\frac{\delta U}{\delta x_{j}} = \rho(x_{j}) \frac{\delta V}{\delta x_{j}}$, 
 that is we can represent $\rho$ by $\rho(x)$ or $\tilde{\rho}(v)$. Assume that $\rho \rightarrow {}^{t} \rho$ 
 preserves type, given $\rho + {}^{t} \rho \rightarrow \infty$ in the $\infty$,
 we have a reduced representation for $U + {}^{t} U$. The condition $X(fg)=0$ means that $f X(g) + g X(f)=0$, 
 that is $-g/f = X(g)/X(f)$. In particular,
 $X(\rho f)=\rho X(f) + f X(\rho)$, that is as long as $X(\rho)=0$, $d U \rightarrow d U^{2}$ preserves type.
 Note that $X(f^{2}) / f^{2} = X(f) / f$

Given $S=\{ f(\zeta)=U f(\zeta) \}$. Consider S for $f \in L_{ac}^{1}$ and continue S to
$\tilde{S}$, so that $\tilde{S} \rightarrow \tilde{U}$ continuous. Given the continuation algebraic,
we have that $\tilde{U} \subset U$ implies $S \subset \tilde{S}$. Given $R(\tilde{U})$ contractible, 
the continuation is simply connected, that is an approximation property with preservation of type.
 
 Assume that $\frac{d \Phi}{d t}=d U(\Phi)$, we then have $\frac{d \Phi}{d t}=\xi(t) \frac{\delta \Phi}{\delta x} + \eta(t) \frac{\delta \Phi}{\delta y}$,
 further $\int \frac{d \Phi}{d t}=\int \Phi d U \simeq U \Phi$, where U is taken as a functional and 
 given $\Phi$ absolute continuous, we have $\Omega(d U)=\Omega(U)$.
 Any constant surface to a holomorphic function, corresponds to a multivalent surface.
 A Pfaff condition $d U = \rho d V$ maps planar regions onto planar regions.
 Given $U_{j}^{*} \Psi_{k}$ algebraic ($U^{*}$ refers to ${}^{t} U$ with respect to the algebraic dual to $d U$), we have $\int_{\Omega} \Psi_{k} d U_{j}=0$, implies 
 $m \Omega=0$, that is Lebesgue measure zero. 
 
 In particular, in the case $d U (\Psi_{j}) = d UW =0$, with $\Phi \Psi_{1} = \Psi_{2}$, for instance 
 $\Phi(d W_{1})=\rho d W_{1}=d W_{2}$, we have that over $N(d UW)$, we must have $\Phi=C$. 
 Note that $N(d UW) \subset \Omega(\Psi_{k})$.
 Further, note that given $\Phi$ algebraic, it must map constants onto constants.  In a two-mirror model 
 ($\Psi_{1} \rightarrow \Psi_{2} \rightarrow \Psi_{3}$), 
 we must have $\Phi_{1} \Psi_{1} = \Phi_{3}^{-1} \Psi_{3}$.
 For $N(d UW)$, we must have $N(d UW) \subset \cap \Omega(\Psi_{k})$, that is not all $\Phi_{k}$ 
 must be invertible.
 
 Example: Assume for instance, that every f in R(UW) is such that $UWf=UW_{j} f$, for some j. 
 Assume that $U_{j} \bot U_{k}$ and $U_{j}^{2}=U_{j}$ over R(W), 
 then we can determine f uniquely, where $U_{j} W_{j} f=f$.
  
 Assume that $U f=\int f d U= \int f d I=If$ 
 $\sim \int U d f = \int I df$. Given $U_{j}$ analytic over $df \neq 0$, with $ U_{j} d f=0$ and $U_{j} \rightarrow U$, assuming uniformity,
 we have that $(U-I) d f \equiv 0$ on $\Omega$ or $m \Omega =0$ (algebraic). 
  
  \subsection{Relative invariance}
 Assume (\cite{Lie88}) $\forall x,y \in D$ a domain, there is $\Phi \in \mathcal{G}_{r}$, with $\Phi(x)=y$, then
 $\mathcal{G}_{r}$ is called transitive.
 For a transitive group $\mathcal{G}_{r}$ in the n-space, we have $r \geq n$. Given r=n
 the group is simply transitive, we have in this case $d I \in \mathcal{G}_{r}$.
 
  $T_{b}$ is invariant for $T_{a}$, if $T_{a}^{-1} T_{b} T_{a} = T_{c}$, where $c=c(a,b)$ (\cite{Lie88}, kapitel 15).  
  Note that normal operators are consequently relative invariants and in $\mathcal{G}_{1}$ exchangeable.
  Given two groups $\mathcal{G}_{\mu},\mathcal{G}_{m}$, invariant sub-groups of $\mathcal{G}_{r}$,
  with $\mathcal{G}_{\mu} \cap \mathcal{G}_{m}=\{ 0 \}$, then $\big[ X_{\mu}, X_{m} \big]=0$ (\cite{Lie88}, kapitel 15).
 
  Coefficients to $d U + d U^{\diamondsuit}$ are given by $((\xi - \eta),(\xi + \eta))$, with involution condition
 $\frac{\xi - \eta}{\xi + \eta} = -\frac{d x}{d y}$. Note in particular $(\xi - \eta)^{2} + (\xi + \eta)^{2} = 
 2 (\xi^{2} + \eta^{2})$, that is standard complexification. Consider the conjugation 
 $d U^{\bot} = \sigma d U^{\diamondsuit}$ with $\sigma \sim 1$. Given $d U^{-1 \diamondsuit} U^{\bot} U^{\diamondsuit}=
  \alpha d U^{\bot} U^{\diamondsuit}=
  \alpha \beta d U^{\diamondsuit}$, where $\alpha \sim 1 / \beta$, we have $U^{\bot} \sim U^{\diamondsuit}$.
  Given $d U,d U^{\diamondsuit}$ both closed, we can consider d U as harmonic.
  
 A perfect transformation dV is given by $\big[ dV, dU_{k} \big]=0$ $\forall k$, For instance dV 
 regular and perfect, implies that the base
 for  $\mathcal{G}$ has a regular representation. Given $\mathcal{G}_{r}$ is primitive, there are no perfect transformations.  
 Consider partial regularity, for instance $\mathcal{G}_{r} \subset \mathcal{G}_{q}$ with $\big[ dV, dU_{k} \big]=0$ 
 for $k=r+1,\ldots,q$. Note that given $\big[ dV ,dU_{k} \big]=0$, with coefficients $\xi,\xi_{k}$, 
 we have $dV(\xi_{k})=dU_{k}(\xi)$,$dV(\eta_{k})=dU_{k}(\eta)$.
 Given d V perfect with respect to d U, we then have $d V$ absolute continuous relative d U and 
 conversely. 
 
  Given U projective with $U^{2}(a) = U(b)$, we have $\Omega(a) \simeq \Omega(b)$.
  U is projective relative R(V), if $U^{2}V=UV$, that is UV=V.
  Note that $d U^{\bot}(U)$ of bounded variation, does not imply $d U$ of bounded variation, for instance 
  $\frac{d U^{\bot}}{d U} \frac{d U}{d T}=\rho \frac{d U}{d T}$.

  \subsection{Invertibility}

 Consider $(x,y)$ in a simply connected domain.
 Assume that $\rho d U^{-1}(x,y)=dU(1/x,1/y)$, or for $f=e^{\phi}$, $\widehat{d U}(-\phi)=\rho \widehat{d U^{-1}}(\phi)$.
 Given $d U(1/x,1/y) \rightarrow C$ regularly for finite arguments, we have when $\rho$ constant, that 
 $d U^{-1}(x,y) \rightarrow C$ regularly for finite arguments.
 Assume now $\rho$ regular. Note that in the case where $\rho$ is linear and d U analytic, we have existence of 
 $dU^{-1}$ analytic according to the equation. More precisely, assume that $d U = \rho d V$ in $\mathcal{E}^{'(0)}$, in analogy with Poincare, assume that $\rho$ analytic with absence of 
 essential singularities in $\infty$, we can then find d U and d V analytic, such that $d U=\rho d V$. 
 When $\rho$ can be chosen as analytic, the equation preserves 
 harmonicity. For instance, assume
 $d U(f) = \rho d V(f)$, with $\rho f \in C^{1}$. Given $d V > 0$ and $d V \rightarrow d I$ regularly, we have 
 existence of 
 $\lim_{\rho > 1} U = \lim_{\rho < 1} U$ and we have $UI=IU \in G$. The condition implies invertibility (orientability).
 Consider a domain for analyticity $\Omega$, for both $d U,d V$ simultaneously. Continue this with 
 complex lines to the infinity, that is $\tilde{\Omega}$ simply connected. Given $R(V)$ planar  and dense 
 in the domain for $d U$, it is sufficient to
 consider $\rho$ linear. Given d U-d V of bounded variation on such a line L, in a planar domain, we must have d U=d V on a disk,
 and when L can be chosen arbitrarily, $d U \simeq dV$. Conversely, given $\rho=1$ has finite order, 
 then d U,d V are not necessarily independent
 (or dependent). The condition $\int_{\Omega} \rho d V=0$ implies $m \Omega=0$ thus means on 
 $\{ \rho = 1 \}$, that $d U \neq 0$ and $d V=0$ implies $m \Omega=0$. 
 
 Note that the cylinder web can not be generated by translation, thus the convex closure of a figure,
 must be defined relative the group G. We refer to a linearly convex continuation, when the movement
 preserve disks. Consider the convex closure according to $\mid U_{1} -I \mid \leq 1,
 \mid U_{2} - I \mid \leq 1$, then obviously the cylinder with web $U_{1}=U_{2}$ is contained in the
 linear convex closure T(C). Note that when U is one-parameter, it does not imply that
 ${}^{t} U$ is one-parameter, if not $U \rightarrow {}^{t} U$ preserves type. Note that $\Sigma=\{ d F=0 \}$ 
 and $\tilde{\Sigma}=\{ d^{2} F=0 \}$, we have that
 given a figure S, such that $S \cap \Sigma=\emptyset$, we may still have that $S \cap \tilde{\Sigma} \neq \emptyset$
 in particular, movements that are solvable on S, do not generate S.

  \subsection{Symmetry}
Assume that $d H=-\xi Y + \eta X$.
Assume that $X,Y$ are analytic, we then have $\Delta H(f)=- (\overline{d} \xi) Y + (\overline{d} \eta) X$.
Immediately, $<d H(f), \phi>=<f, (\xi_{x} + \eta_{y}) \phi> + <f, d H(\phi)>$, that is a condition for $d H \rightarrow d {}^{t} H$
to preserve type, is that $(\xi_{x} + \eta_{y})=0$, given $(\xi,\eta) \rightarrow {}^{t} (\xi,\eta)$ preserves type.

Consider $<d U^{\diamondsuit} f, \phi>=<f, (-\eta_{x} + \xi_{y}) \phi>
+ <f, d U^{\diamondsuit} \phi>$, that is the orientation can be determined by $d U \rightarrow d U^{\diamondsuit}$.
If we consider a minimal operator, $d U^{\diamondsuit} \phi$ and a maximal operator $d U^{\diamondsuit} f$, 
we see that realizations are dependent on orientation.
In particular $<d U^{\diamondsuit} f,\phi>= <f, d U^{\diamondsuit} \phi> + <d V f, \phi>$ corresponds to 
an oriented realization, where $d V > 0$, as $f > 0$. Given $d U^{\bot} + d V = d U$, we have that 
the condition $d V > 0$ excludes presence of a spiral in R(U). The condition $\int f d U^{\diamondsuit}=0$ can be related
to vanishing flux.  
\emph{Given existence of G with $\frac{\delta G}{\delta x}=-\eta$ and $\frac{\delta G}{\delta y}=\xi$
(for instance $\xi,\eta$ analytic), a condition for orientability is that $\Delta G \geq 0$, 
that is we have a separately absolute continuous Hamiltonian.}
Assume that $\Delta G=0$ on $\Omega (U_{1},U_{2})$ corresponds to an equation $d F=0$ on $\Omega (x,y,z)$. 
Analogously with Green,
$\int \frac{\delta G}{\delta U_{1}} d U_{1} - \int \frac{\delta G}{\delta U_{2}} d U_{2}=0$ iff
$\frac{\delta^{2} G}{\delta U_{1} \delta U_{2}} = \frac{\delta^{2} G}{\delta U_{2} \delta U_{1}}$
and $-\int \frac{\delta G}{\delta U_{2}} d U_{1}- \int \frac{\delta G}{\delta U_{1}} d U_{2}=0$
iff $\Delta G=0$, in particular there is $G_{1} \sim G$ absolute continuous in $(U_{1},U_{2})$ 
separately. 

Example: consider translation $\tau$ and $f(\tau x, \tau y) = f(x,y) + C_{x} + C_{y}$, that is $C_{x}=-C_{y}$ gives the corresponding invariants.
Assume that $(d f)(\tau x,\tau y) = (d f)(x,y)$ that is $C'_{x}=-C'_{y}$, which is the case when $f \in L^{1}$
with a separations axiom of Schwartz type.
Note that $d f^{2} \sim 2 f d f$, that is $f^{2}=const$ does not imply $f=const$. 

Consider $d U^{\bot} = \rho d U$ with $\Gamma=\{ d U = d U^{\bot} \}$ and $\Sigma=\{ d U^{\bot} - d U > 0 \}$ 
a Riemann surface, that is we assume $\frac{d}{d T}(U^{\bot} - U) \neq 0$
regularly, for instance given $\rho$ sub harmonic and locally 1-1 on $\Sigma$ parabolic, 
we have that $d U^{\bot} \sim d U + d I$, that is d U is ``projective'' over $\Sigma$. On $\Gamma$ d U is no longer algebraic, 
that is $\overline{\Sigma}=\Sigma \cup \Gamma$, does not have the approximation property . 
Sufficient for an algebraic boundary, is $d U^{\bot} I - d U I = d I U - d U^{\bot} I$, for instance 
$d U^{\bot} I = - d I U$, $d U I = - d I U^{\bot} $.

Example: assume p a polynomial and $T=pI$,${}^{t} T={}^{t} p I$, then T is a normal operator, with
$IT(\phi)=<T(\phi), d I>=<dI, T>(\phi)=TI(\phi)$.

\section{Classification of groups}

\subsection{Subgroups}

Note that $\mathcal{G}_{w_{1}} \subset \mathcal{G}$ implies for invariant sets that $\Omega \subset \Omega_{w_{1}}$. 
Given $\mathcal{G}_{r}$ in an invariant point P, with $w_{1}$ independent transformations, we have 
 existence of $\mathcal{G}_{w_{1}} \subset \mathcal{G}_{r}$ (\cite{Lie88}, kapitel 12). If $\Omega_{j}$ are invariants sets for $X_{j} \in \mathcal{G}_{r}$,
 independent for a m-group, we have $\mathcal{G}_{m} \subset \mathcal{G}_{r}$. Given two sub-groups
 $\mathcal{G}_{\mu},\mathcal{G}_{m}$ to $\mathcal{G}_{r}$, and if we have existence of $m+\mu-r$ independent 
 transformations in $\mathcal{G}_{\mu} \cap \mathcal{G}_{m}$, they generate a sub-group (\cite{Lie88}, kapitel 12).
 Note that $\mathcal{G}_{w_{1}}(u_{1},u_{2}) \subset \mathcal{G}_{r}(u_{1},u_{2})$, which does not imply 
 $\mathcal{G}_{w_{1}+2}(x,y,z) \subset \mathcal{G}_{r+2}(x,y,z)$.
 Example: the cylinder web in (x,y,z) does not contain the origin, that is the spiral does not correspond to
 a sub-group, but through change of coordinates $(x,y,z) \rightarrow (u_{1},u_{2})$, spirals $u_{1}=u_{2}$ 
 can be considered as a sub-group.

 Type of movement is determined by invariants $\Omega$ over a domain where 
 the measure is analytic. A measure represented in $H'$, is considered on convex neighborhoods 
 of the invariants sets.
 Spiral example: $d U^{\bot}=\rho d U$, with $\mid \rho-1 \mid \leq 1$, where we assume $G \simeq G^{\bot}$,
 has a projective correspondence $\Sigma \rightarrow \Sigma^{\bot}$. 
 
 Example (n=3): translation does not change character when the parameter varies. $\mathcal{G}_{1}$ is 
 of order 1 and given f defined on a planar domain, there is no non-trivial sub-group. Note that
 a continuation $\mathcal{L}_{ac} \rightarrow \mathcal{L}_{c}$
 can be dependent of $u_{3}$ (scaling). 
 
 Example: Assume that $f_{j}$ are m+1 eigen vectors, and assume existence of 
 $\Phi_{j} : f_{j} \rightarrow f_{j+1}$, through a first trivial mapping, thus we have m transformations
 that form a m-mirror model. 
 Further, $U \Phi_{1} f_{1} \sim \Phi_{1} f_{1}$ and $ U  f_{1} \sim f_{1}$, that is 
 $U \Phi_{1} = \Phi_{1} U$ and $\Phi$ is invariant for U and vice versa. Given eigen vectors 
 defined on planar sets through 0 (independent of scaling), we have a m-mirror model.

 According to Lie (Sats 1: \cite{Lie88}) $d U(f)=f_{i}(u,a)$ is of order r iff we do not have existence of 
 $\chi_{k} \neq 0$, such that $\Sigma_{k=1,\ldots,r} \chi_{k}(a) \frac{\delta f_{i}}{\delta a_{k}}=0$ 
 $\forall i$, that is the proposition in the theorem is that
 given maximal order for movement, there are no invariant sets. Note however, if $d WV=\rho d V$
 and if $\rho=const$ gives an extension of the domain, for instance $\rho = 0$ on $\Gamma$, 
 where $d V \neq 0$ on $\Gamma$, then d WV does not have maximal order on the extended domain.

Note $\int f d U =\int \Sigma \frac{\delta}{\delta a_{k}} (Uf) d a_{k}$. Consider now movements that 
can be given as reflection of one-parameter movements, that is $\xi_{1} / \eta_{1} \sim \xi_{2} / \eta_{2}$. 
Given $d U_{k}$ a complete system $k=1,\ldots,q$, we then have $d V f=\Sigma \eta_{i}(x) \frac{\delta f}{\delta x_{i}}$
one-parameter iff $\big[ d U_{k}, d V \big]=0$, that is characteristics to $d U_{k}$ generate invariants to 
d V. Thus, as long as we are considering one-parameter movements, the type is given by characteristics to $d U_{k}$.

\subsection{Convexity}

Given uniform convergence, we can assume that any singularity can be given on first surfaces.
Note that $\int (U^{*} \xi_{j})' d x=\int \xi_{j} d U$, that is for a movement, 
to preserve type for coefficients, $\int \xi_{j} d U=U^{*} \xi_{j} \simeq \xi_{j}$, $j=1,2$,
$U \rightarrow U^{*}$ must preserve type. 

Given $R(d H)$ has the approximation property , consider $\{ Hf < \lambda \} \subset \subset \Omega$
(an algebraic polyhedron), we then have existence of an analytic lifting transformation $d H \rightarrow d U$.
More precisely, relative $\Delta U=0$ or given $d U = \rho d H$, relative $\rho$ analytic.
Consider $\{ \rho=1 \}$
as a hyper plane, then the range for $R(d U)$ can be divided into $\rho  > 1$ and $\rho < 1$. 
Assume that $I \sim \frac{1}{2}(I_{-} + I_{+})$,
that is (linear) convexity implies presence of at two-sided continuous limit. The condition
$\frac{\delta \rho}{d T} > 0$ implies U convex with respect to H, that is given $\rho$ monotonous, 
the range can be given relative a hyper plane. Note further that as $Hf$ is defined by $\int f d H$, if
$d U = \sigma d I$, with $\sigma \rightarrow 1$ regularly, then $H f \rightarrow f$, when $\sigma / \rho \rightarrow 1$. 
In particular, an approximation property for $R(d H)$ implies an approximation property for $R(H)$.

\emph{
A separating functional to $\rho(u_{1},u_{2})$ does not imply a separating functional in $(x,y,z)$.
A sufficient condition for preservation of connected sets $(u_{1},u_{2})$ to 
$(x,y,z)$, is that $<d U,f> \in H(x,y,z)$ (that is analytic).}

Assume that $d U / d V=\rho \geq 0$ with $\rho \in L^{1}_{ac}$, that is U monotonous relative V.
Note that where $\rho'(v) \geq 0$, U can be seen as convex relative V. Where $\overline{\delta} \rho (x,y) =0$,
the continuation preserves harmonicity.

 Assume that $d U_{r}=\rho d H$, that is $H=U_{w_{1}} \in G_{w_{1}}$,
 with $\rho > 0$. We have $\mathcal{G}_{w_{1}}$ represents a hyper plane through 0.
 The condition $\rho' > 0$ implies that the continuation preserves orientation.
 Example : $d U(x,y,z)=
 \rho(x,y,z) d I(x,y,z)$
 becomes $d U(h_{j})=$ $\rho(h_{j}) dI (h_{j})$ and 
 $\frac{\delta \rho}{\delta x}=\Sigma \frac{\delta \rho}{\delta h_{j}} \frac{\delta h_{j}}{\delta x}$,
 with $\rho$ constant $\forall h_{j}$, implies $\rho$ constant in x,y,z.

Given $d U F=\xi \frac{\delta F}{\delta x} + \eta \frac{\delta F}{\delta y}$ and $d A G=X \frac{\delta G}{\delta x} + Y \frac{\delta G}{\delta y}$.
We then have $d U(d AF) - d A(d UF)= (d UX - d A \xi) \frac{\delta F}{\delta x} + (d UY - d A \eta) \frac{\delta F}{\delta y} + R$
where $R= \xi Y \frac{\delta^{2} F}{\delta x \delta y} + \eta X \frac{\delta^{2} F}{\delta y \delta x}
- \eta X \frac{\delta^{2} F}{\delta x \delta y} - \xi Y \frac{\delta^{2} F}{\delta y \delta x}$. Thus
given $F \in C^{2}$ and $\frac{\delta^{2} F}{\delta x \delta y} = \frac{\delta^{2} F}{\delta y \delta x}$ 
we have R=0. 
Given $d U(d A) - d A (d U) = \lambda d A$ defines a one-parameter group, assume that dU preserves type 
for X,Y and dA preserves type for $\xi,\eta$, then the condition is given by $d A \xi / d A \eta \simeq \xi / \eta \simeq X / Y$.

Assume that $\big[ d U_{i},d U_{j} \big](f)=(d U_{i}(\xi_{j}) - d U_{j}(\xi_{i})) \frac{\delta f}{\delta x} + 
(d U_{i}(\eta_{j}) - d U_{j}(\eta_{i})) \frac{\delta f}{\delta y}$. Analogous, given $d U_{i} = \rho_{i}d I$,
we have $\frac{\delta g}{\delta y}=\rho_{i} \xi_{j} - \rho_{j} \xi_{i}$ and $- \frac{\delta g}{\delta x}=\rho_{j} \eta_{i} - \rho_{i} \eta_{j}$.
Assume further that $c_{1} \frac{\delta g}{\delta y}=\{ h,v \} \frac{\delta f}{\delta x}$ and $c_{2} \frac{\delta g}{\delta x}=\{ h,w \} \frac{\delta f}{\delta y}$,
for constants $c_{1},c_{2}$, that is $\{ h,v \}$ is related to $\{ h,w \}$ as $d x/ d y$ is related to $d y / d x$.
Given that translation is analytic over g, we have $\big[ d U_{i}, d U_{j} \big] \simeq 0$. Conversely, given
$\big[ d U_{i}, d U_{j} \big]=0$, we have $\frac{\{ h,v \}}{ \{ h,w \}} \simeq \frac{d x}{d y}$.
Given h,v functions, $<\{ h,v \},\phi>= <\frac{\delta v}{\delta y},\frac{\delta}{\delta x}(h \phi)> - 
<\frac{\delta h}{\delta y}, \frac{\delta}{\delta x} (v \phi)> - <\frac{\delta v}{\delta y}, h \frac{\delta \phi}{\delta x}> + 
<\frac{\delta h}{\delta y}, v \frac{\delta \phi}{\delta x}>$ or $=<\frac{\delta v}{\delta y},\frac{\delta }{\delta x}(h \phi)> -
<\frac{\delta v}{\delta x}, \frac{\delta}{\delta y}(h \phi)> +  <\frac{\delta v}{\delta y} , h \frac{\delta \phi}{\delta y}>$ -
$<\frac{\delta v}{\delta y}, h \frac{\delta \phi}{\delta x}>$  
, that is the equation =0, if we have 
symmetry for $h \phi$ and $\frac{\delta \phi}{\delta x}=\frac{\delta \phi}{\delta y}$, that is $\phi$ symmetric.

Note that given $Uf$ absolute continuous, we have $U^{2}(f)=\int \frac{d U^{2}(f)}{d T} d T \sim$ 
$2 \int (U f) \frac{d U}{d T} d T=$ $2 \int \int d U (d U(f))$.
Assume that $\frac{d c}{d a}=\frac{d c d b}{d b d a}$, with the condition that $d c/d a=1$, implies $c \in A$. 
 Note that c is absolute continuous in a does not imply that b is absolute continuous in a. Given d c/ d b = d a / d b 
 then we must have $c \in A$. Further, we can determine $\theta$ so that
d U is absolute continuous in $\theta$, but also d U is absolute continuous in x,y,z separately. Given F 
is absolute continuous with respect to $\theta$, $F=\int \frac{d F}{d \theta} \frac{d \theta}{d u} d u$,
where $\frac{d \theta}{d u} d u$ of bounded variation and monotonous, that is according to the Stieltjes integral.

Consider $P (x,y,z) \rightarrow P(u_{1},u_{2})$. Given $\frac{d U}{d U_{1}}=\rho$ and $\frac{d U}{d U_{2}}=\vartheta$
with $\rho,\vartheta \in L^{1}$, we can determine a maximal domain for $\rho,\vartheta$ absolute continuous. 
Assume $d U \rightarrow d I$ regularly and that we have existence of $d U=\lambda d V$ and $d V=\sigma d I$.
Then convexity requires a two-sided limit for $dV$. For instance, if $\lambda \rho > 0$, then $\lambda,\sigma$
must have the same sign and we have two possibilities for sign of $\lambda$ (and $\sigma$). 
When $\rho$ is absolute continuous and $\geq 0$, we have $d \rho=0$ implies $d U$ convex iff $d U_{1}$ 
convex.

  \subsection{Orientation of movement}

 Assume that P is an invariant point to $\mathcal{G}_{m},\mathcal{G}_{r}$, with $\mathcal{G}_{m} \cap \mathcal{G}_{r} = \{ 0 \}$, 
 we then have $\mathcal{G}_{m} \simeq \mathcal{G}_{r}$
 Assume that $\mathcal{G}_{m} \cap \mathcal{G}_{r} \neq \{ 0 \}$, why
 $\big[ X_{m}, X_{r} \big]$ defines a new movement Y. Assume that $\big[ Y, X_{m} \big] \neq 0$
 and $\big[ Y, X_{r} \big] \neq 0$ According to Jacobi: $0=\big[ Y, \big[ X_{m}, X_{r} \big] \big] + 
 \big[ X_{m}, \big[ X_{r},Y \big] \big] + \big[ X_{r}, \big[ Y,X_{m} \big] \big]$. Thus
 $\big[[ Y, \big[ X_{m},X_{r} \big] \big]=0$ implies $\big[ Y,X_{m} \big] \sim \big[ Y, X_{r} \big]$.
  
  Assume that $a^{0}$ fixed, it generates a movement $f_{i}(x,a^{0})$. Assume that the movement 
  preserves character in a neighborhood of $a^{0}$ (disk). Preservation of type
  corresponds to an involution condition $\xi/\eta \simeq - \frac{d x}{d y}$. Consider $d U / d T$
  and $U \rightarrow (x,y)_{T}$, this defines a neighborhood of $(x,y)$, as T varies. 
  More precisely, given $d U = \rho d V$, it is sufficient for preservation of regularity, that
  $\rho,1/\rho$ both are regular, that is on the domain for simultaneous analyticity, we have that 
  $d U =0$ iff $d V=0$. Assume that $\mathcal{G}_{r}$ has a local pseudo base over $f=e^{\phi}$ with 
  $\phi \mid_{\Omega_{j}} = \phi_{j}$ patch functions and $\phi_{i}=\phi_{j}$ on $\Omega_{i} \cap \Omega_{j}$. 
  Sufficient for preservation of orientation, is for instance
  $\eta_{x} - \xi_{y} > 0$, that is $\eta_{i} / \xi_{i} \simeq \eta_{k} / \xi_{k}$ on $\Omega_{i} \cap \Omega_{j}$
  further for instance $\frac{(\eta_{i})_{x}}{(\xi_{i})_{y}} \simeq \frac{(\eta_{k})_{x}}{(\xi_{k})_{y}}$.
  
  According to Lie (Theorem 6: \cite{Lie88}), given $d U=f(u,a)$ one-parameter, so that $d I \in \mathcal{G}$, 
  the transformations can be ordered as pairwise inverses. Further we have existence of a unique 
  transformation, dependent of parameter, that generates $\mathcal{G}$.
  According to kapitel 4 (\cite{Lie88}), given a r-parameter group including the id-transformation, 
  there is a base of independent transformations in one parameter. Note that when $d U=\rho d U_{1}$,
  where $(d U_{1},d U_{2})$ is dense in the domain for $d U$, where $d U$ is independent on $d U_{2}$,
  the orientation for d U, can be considered as single valued.

 \subsection{A weak type of movement}
 
 Assume that H is a hyper plane through 0, we discuss movements Y symmetric with respect to H.
 Since $\big[ X,X \big]=0$, we can given $\big[ X,Y \big]=H$, assume that $\big[ X,H \big]=0$, that is
H can be given as reflection through X or Y.
 For instance $\big[ H,Y \big]=- \big[ Y , H \big]=0$.
 Given $Y=H \pm X$, we then have $\big[ X,H \big] = -\big[ H,X \big]$. Example: pure transformations:
 given that conjugation $X \rightarrow X^{*}$, is related to harmonic conjugation, $X^{*} \sim X^{\diamondsuit}$,
 pure forms $(X,X^{*}) \rightarrow -(-X^{*},X)$, are related to uniformity.

Invariant factors (\cite{Riesz37}) 
We consider convex bodies $d U = \lambda_{1} d U_{1} + \lambda_{2} d U_{2}$, with $\lambda_{1} + \lambda_{2}=1$, 
$\lambda_{j} \geq 0$. Consider modules $M_{j}$ in n-space. We write $(b) \subseteq (a)$, if $a \mid b$. 
Assume that $(e_{k})=\{ x \in M \quad L_{v_{j}}=0 \quad v=v_{j} \quad j=1,\ldots,k \}$, where $v_{j}$ are integers and 
$L_{j}$ are linear functionals. Modules $(e_{j})$ are formed as linearly independent. $(e_{k})' = \{ L_{v}(x) \quad x \in (e_{k}) \}$. 
Then we have existence a minimal module $(e_{j}) \subset (e_{k})'$.
Thus, when $d U \rightarrow d I$ regularly, $d U_{2}$ can be constructed with a maximal domain for convexity and relative a scalar 
product, for instance as $d U_{2} = d U_{1}^{\bot}(\lambda_{2})$.

A decomposition lemma: (\cite{Lie88} kapitel 24) Assume that $M_{r}$ invariant for $\mathcal{G}_{r}$, 
we then have $M_{r} \subset M_{r-s+h} \subset M_{r-s}$ and
$\mathcal{G}_{r-s} \subset \mathcal{G}_{r}$. Choose $P \in M_{r-s} \backslash M_{r}$ and $\Phi$ a transformation 
to $\mathcal{G}_{r-s+h}$, this gives a h-fold image of $\Phi P$. Given $\Psi$ a transformation to 
$\mathcal{G}_{r}$, this gives a s-h fold image of $\Psi M$. 
This gives all for $\mathcal{G}_{r}$ invariant decompositions, when we let $\mathcal{G}_{r-s+h}$ 
vary. 
 
In particular, $d V \mid d U$ implies $(d U) \subset (d V)$ in the sense of sets. Assume that 
$d U=\rho d V$, that is given $dU(f)=0$ and $\rho \in L^{1}(d I)$, we have existence of $d V \in \mathcal{G}$.
Where $\rho \neq 0$, we have existence av $1/\rho$, that is $d V=\frac{1}{\rho} d U$. Compare with 
Nullstellensatz, that is given $\mid d U(f) \mid \leq C$ and $\rho$ downward bounded, then 
$\mid d V \mid \leq C'$. A sufficient condition for this, is that $\{ \rho < \lambda \} \subset \subset \Omega$
on the domain for d V. 
Note that given $\rho \rightarrow 0$ and $d U \bot d U^{\bot}$, we do not have $d V \bot d V^{\bot}$.

Consider $d U=\rho_{1} d U_{1} + \rho_{2} d U_{2}$, with $\rho_{1} + \rho_{2}=I$, that on 
first surfaces to $\rho_{j}$, give convex bodies. Assume that $U=\rho_{1} U_{1} + \rho_{2} U_{2}$, 
we then have $d U=\rho_{1} d U_{1} + \rho_{2} d U_{2}$ over
$\{ \rho_{1} =const,  \rho_{2}=const \}$ and given $\rho_{j} \rightarrow \rho_{j}^{*}$ preserve type 
also the converse holds. More precisely, consider $d U = \rho_{1} d U_{1} + \rho_{2} d U_{2}$, with 
$\rho_{1} + \rho_{2}=I$ as conjugation. Given $U=U(U_{1},U_{2})$, we have $\frac{d^{2} U}{d T^{2}}=(\frac{d^{2} U}{d U_{1}^{2}}) (\frac{d U_{1}}{d T})^{2} +
(\frac{d^{2} U}{d U_{2}^{2}}) (\frac{d U_{2}}{d T})^{2} + \rho_{1} \frac{d^{2} U_{1}}{d T^{2}} + \rho_{2} \frac{d^{2} U_{2}}{d T^{2}}$.
As $\rho_{j}$ constants, only the two last terms remain.

\newtheorem{Main}[prel]{Proposition}
\begin{Main}
A weak type of movement can be defined by invariant factors on a non-trivial domain, 
through $<d \lambda_{1} U_{1} + d \lambda_{2} U_{2},\phi>=\lambda_{1} <d U_{1},\phi> + \lambda_{2} <d U_{2},\phi>$
$=<d I,\phi>$. Assume either of the movements is analytic over $\phi$,
then the weak type coincides with the type, only when $\phi$ is absolute continuous. When $d U \rightarrow d {}^{t} U$
preserves type, the argument can be continued to $T(\phi)$ for $T \in \mathcal{H}(E)$
\end{Main}

Example : $UV e^{\phi}=e^{(U + V) \phi}$ with 
U maximal, absolute continuous over V. 
Given projectivity, that is U+V=I, we have UV=I, that is where (U,V) is dense in the domain for $e^{\phi}$, 
then U is invertible. Where U is absolute continuous relative V, we have $\frac{d U}{d V}=0$ implies $UV=V$, 
that is U is a continuation of $R(V)$.

Example: Assume that $<d V_{1},d V_{1}>=1$ and $<d V_{1},d V_{2}>=0$, and $d U=\rho_{1} d U_{1} + \rho_{2} d U_{2}$,
Assume further $<d V_{1}, d V_{2}>=- <d V_{2},d V_{1}>$ and $<\rho_{1} d V_{1}, d V_{1}>=\vartheta_{1}$,
we then have $<d U, d U>=\vartheta_{1}^{2} + \vartheta_{2}^{2}$ ($=1$).
Assume further that $\overline{\vartheta_{j}} =<d V_{j}, \rho_{j} d V_{j}>$, with compatibility conditions
$\vartheta_{1} \overline{\vartheta_{2}} - \vartheta_{2} \overline{\vartheta_{1}}=0$, we get convexity
for the movement, a contact transform according to Legendre.
 
 Assume that $\mathcal{G}_{r}$, with $n=m+r$ where $\Omega_{m}$ gives continuation of $\mathcal{G}_{r}$ 
 to the n-space, for instance $\mathcal{G}_{r} \ni d U$ with $d U=d I$ over $\Omega_{m}$, or a zero-space 
 to the movement in the phase. Given that U=I over a convex set with respect to $\mathcal{G}_{m}$, 
 gives existence of $d W \in \mathcal{G}_{m}$ regular, with $d UW=0$ implies $UW=W$, that is invariants 
 form a normal tube. When $R(W)=\cup L_{j}$, where $L_{j}$ is $R(G_{j})$ one-parameter, we have existence
 of approximation with single valued orientation.

 \section{Compact sets}
 
 Assume that $E_{n}$ euclidian and finite dimensional. Assume that $E_{R}=\{ x \notin E_{n} \quad \mbox{ dist }(x,E_{n}) > R \}$.
Lebesgue's theorem gives that, when $\parallel f \parallel \leq M$, we have uniformly, that $\parallel f(x+h) - f(x) \parallel \rightarrow 0$, 
as $\mid h \mid \rightarrow 0$ and $\parallel f \parallel_{E_{R}} \rightarrow 0$, as $R \rightarrow \infty$,
 iff $\{ f \} \subset L^{p}$ is $L^{p}-$ compact. Hausdorff compactness is equivalent with $L^{p}$-compactness
(\cite{Riesz33})

 \subsection{Compactness relative sub groups}
 
Consider $F(u_{1},u_{2})$ as a planar function, corresponding to pseudo-vectors, that is with double orientation.
In $L_{c}^{1}$ we consider $F(u_{1},u_{2},u_{3})$, where $u_{3}$ refers to scaling. For absolute continuous measures 
we can use planar functions.
Note for closed forms, we have independence of $u_{j}$ (independence of local 
coordinates).

Consider the pseudo inverse: $(U,+)^{-1} \simeq (U,-)$, where we assume $\{ d (U,+) < \lambda \} \subset \subset E_{n}$
if $1/ d (U,+) \rightarrow 0$ i $\infty$, in particular this condition means absence of essential 
singularities in the $\infty$.
Example: $\frac{d}{d t} \log (U,+)=\frac{d (U,+)}{(U,+)} = (U,-)d (U,+)$ and $\frac{d}{d t} d (U,+)^{2} \sim 2 (U,+) d (U,+)$.

More precisely, assume that $d {}^{t} U = \vartheta d U$
and $d U = \rho d I$, assume that $\vartheta$ analytic, then according to Poincare, there are $\rho,\phi$ so that
$\vartheta = \phi / \rho$. Starting from
$\vartheta$ analytic, we can thus find $d U,d {}^{t} U$ analytic.  
Assume that f has one-sided support, we do not necessarily have $U_{+} d U_{+} f \simeq U_{-} d U_{+} f$.
Note that given $\Omega_{+}=\{ U_{+} f=f(u) = f(1) \}$ and $\Omega_{-}=\{ U_{-} f= f(1/u) = f(1) \}$,
then we have a two sided limit, when $\Omega_{+} \simeq \Omega_{-}$.
Given $\vartheta$ linear, we can determine $\rho,\phi$ using singular integrals. In the continuation
$\mathcal{L}_{ac} \rightarrow \mathcal{L}_{c}$, the dependence is $(u_{1},u_{2}) \rightarrow (u_{1},u_{2},u_{3})$,
and we can not assume $\rho$ linear, but absence of essential singularities gives the result 
as above.  
Given $\vartheta,1/\vartheta$ analytic, we can reverse the transponats. In the case where  
$\vartheta=const$ implies $(x,y,z)=P$ a point in finite space, we write $d {}^{t} U \bot d U$ and the invariant sets are considered as having a discrete
intersection.
When this relation can be continued to $\infty$, the movements are considered as linearly independent.
\newtheorem{kompakt}[prel]{Lemma}
\begin{kompakt}
Assume that $V^{-1} W$ compact in the sense that $m ( x \quad dist_{\Gamma} (V^{-1} Wf) > R) \rightarrow 0$
in the $\infty$, where $\Gamma=\{ \sigma^{-1} \rho=const \}$, that is a characteristic surface, where $d W=\rho d I$ and
$d V=\sigma d I$ and $\sigma^{-1} \rho \in \mathcal{H}$. Given d V analytic, we assume a continuous
dependence of the parameter. 
Thus, when V,W are close, $V^{-1}W$ is close to compact,
For instance $W f= V f + R f$, where we assume R compact, that is V is an approximation of W. 
\end{kompakt}
\subsection{Range of the group}
Note that convex sets of functions $\mathcal{L}_{c}$
are not equivalent with sets of convex functions in $\mathcal{L}_{ac}$.
Given $\mathcal{H}$ quasi-complete (\cite{Schwartz57}), we have $\mathcal{H}_{c}'$ is uniformly convergent on convex, 
compact, equilibrated sets. Assume $\mathcal{H}$ quasi-complete and consider a subspace of $\mathcal{H}$, where we have existence of $\lim_{x \rightarrow a} f(x)$, 
for $f \in \mathcal{H}$ and a finite, then the subspace is quasi-complete. The subspace of $\mathcal{H}$
of functions, that preserve constant value in the $\infty$ (behaves like a polynomial in 0), is 
quasi-complete and we have absence of essential singularities in the $\infty$. 
Compact sets are quasi-complete,
conversely bounded and quasi-complete sets are relatively compact. 

Starting from $<T,d U>(\phi)=<T(\phi),d U>=<T({}^{t} U^{*} \phi),d I>$. Let $\tilde{E}=\{ \phi \quad U \phi \in E \quad \exists U \in G \}$,
where we assume $I \in G$
and that we have a continuous deformation, with $U^{-1} \tilde{E} \subset \tilde{E}$, for instance  $U T \rightarrow T$ in $C^{\infty}$, when $U \rightarrow I$ regularly, for some U. Note that starting from $<T(\phi) , d U>=<T,d U>(\phi)$, 
in order that 
$<T,d U> \in \mathcal{D}_{L^{1}}'$,
with $d U \in \mathcal{E}^{'(0)}$, we must have $T(\phi) \in C^{(0)}$. 

The type of movement is dependent of topology,
that is we start from $\Omega: dU = \rho d I$, with $\rho \in \mathcal{H}$ implies $\rho = const$ on 
$\Omega$. Define for instance $E_{ac}'=\{ dU \quad d U =0 \Rightarrow U=I \}$.
Consider $R : \mathcal{K}_{c} \rightarrow \mathcal{H}_{ac}$ sub-nuclear, given R can be chosen as Id, 
then $\mathcal{H}_{ac}$ is nuclear.
In particular, assume that $\tilde{\mathcal{H}_{ac}}=\{ \phi \in L^{1}_{c} \quad \sigma \phi \in L^{1}_{ac} \quad \exists \sigma \}$,
where $\sigma \rightarrow 1$ regularly.
For a normal operator, we consider $<T(\phi),d U> \rightarrow <d U(\phi), T>=<d U,\phi> (T)$.
Given $d U \rightarrow d {}^{t} U$ preserves type,
we can write $<d U,\phi> \simeq <\phi,d U>$ and
$<<d U,\phi>,T> \simeq <\phi,<T,d U>>=<\phi,<d U,T>>$.

Consider $< U^{-1} F, U \phi>=<F, {}^{t} U^{-1} U \phi>$, that is given ${}^{t} U^{-1} U \sim I$ and $R : U \rightarrow {}^{t} U$,
we have in a weak sense that all movements that preserve type under transponation, have $R \simeq I$. 
Assume that $d {}^{t} U=\vartheta d U$, then the type for ${}^{t} U$ depends on $\vartheta \rho=const$ on $\Omega$
and $\vartheta \rho \in \mathcal{H}$, that is it is sufficient for the type to be preserved, 
that $\vartheta \sim 1 / \rho$ (=const) on
$\Omega$. 

 \subsection{Compactness relative continuations}

 Given analytic leaves for the orthogonal (multivalent), there is a spiral approximation. Given
 $U_{S} \rightarrow I$ over $\phi$, we have $\phi \in C$ (the cylinder web).
 
 $U(T \phi) - T(U \phi)=(U - I)T(\phi) + T(I-U)(\phi)$, given $TI=IT$. Let $\Omega$ be invariant sets to T
 and $\Sigma$ invariant sets to $\phi$, given $\Sigma \rightarrow \Omega$ projective, the type is preserved.
 Consider type in a weak sense, according to $<T_{x},\phi>_{y} \in \mathcal{H}$ or $T * \phi \in \mathcal{H}$.
 Assume that $d U = \rho d U_{1}$, given $\mid \rho \mid \rightarrow 0$, we have $\mid d U \mid \leq C \mid d U_{1} \mid$,
 that is $d U$ absolute continuous relative $d U_{1}$. We can discuss relative compactness, for instance
 $\{ \mid U - V \mid > A \} \rightarrow 0$, as $A \rightarrow \infty$, that is when $V$ is continuous in the parameter, 
 $\{ f \quad V^{-1}U \neq I \}$ $\downarrow 0$, or $V^{-1}U \neq I$ on compact sets. Thus every such 
 continuation of $R(V)$ (non-trivial) is compact.
 Consider $U f = e^{V \phi}$, with $\Omega=\{ Uf = f \}$ and $\Sigma=\{ V \phi = \phi \}$, given $f,\phi$
 in the same topology. When $\Omega$ gives algebraic first surfaces, then $\Sigma$ gives algebraic zero surfaces,
 that is given $V$ algebraic and $\Omega \rightarrow \Sigma$ projective, then U is algebraic. 
 In particular, given V algebraic and Vf + (I-V)f=f,
 when $f=e^{\phi}$, we have $(I-V)f=0$ iff $V \phi = \phi$, that is the conditions implies a projective 
 equivalence at the boundary.

\section{Representation of movements}
Assume that G is an analytic Lie grupp, with $I \in G$. Assume that A an one-parameter subgroup in G, 
and consider (\cite{Garding47}) $T_{A}x=\lim_{s \rightarrow 0} \frac{1}{s} (T(a) -1) x$.
We assume $D(T_{A})$ dense in B, a Banach space.  Assume that $a \rightarrow T(a)$ a bounded operator on 
 a Banach space.
 Given $a \rightarrow b$, we have $T(a) \rightarrow T(b)$. Assume that $T(a)f= f d U$, 
 measures of bounded variation. Assume $M_{A}(f d U)=\lim_{s \rightarrow 0} s^{-1} \int T(a(s)) f$. Given f harmonic in a symmetric neighborhood of $a^{0}$, 
 we have $M_{A}f$ is constant (one-parameter).

 Note, given $T(a)=T(a_{1},\ldots,a_{n})$, we have $(T(a)-I)(x)=x(a_{1},\ldots,a_{n}) - x(0)$, 
 under the condition x(a) analytic and completely dependent of a. Given x absolute continuous, we have
 $x(a) - x(0)=\int_{0}^{a} d x$, where we assume  a in a connected set in the domain.
 G can be generalized to a manifold M in $E^{m}$, in the sense of Whitney (\cite{Garding47}). 
 
 Assume that $T(a)$ is defined modulo $C^{\infty}$. Continue $T(a)$ to $\mathcal{L}_{c}$, 
 for instance $T(a)$ very regular, as $s \rightarrow 0$. We have $T_{A} f \in C^{\infty}$ over $f \in C^{\infty}$ 
 given $T_{A}$ has algebraic symbol. In this case, given $\Omega=\{ \zeta \quad T(a)f=f \}$ and $\Omega_{s}=\{ \zeta \quad T_{A}f=0 \}$, 
 we have $\Omega_{s} \simeq \Omega$, that is we have independence of s.

 Consider $Uf-If(x)=f(u)-f(1)$. Assume that $U-I=V$ with $d V \in G$, this means that 
 $f(u)-f(1)=0$ on a line L iff $f(v)=0$ on L, that is L is a zero-line to a movement V. Assume that H 
 corresponds to a line through 0 (that is L), we then define $u \rightarrow v$ as conjugation.
 Consider starting from $If = \int f d I$, a $L^{p}$ compact neighborhood (or domain of holomorphy), 
 where we assume for instance $u > h \geq 0$. Assume vanishing flux according to $\int_{L} f(v) d v=0$. 
 Note that when H is considered modulo regularizing action, $H(\phi) \in \mathcal{D}^{'F}$, that is
 of finite order. Further, $\mathcal{D}_{c}^{' m}$, has the strict approximation property and an approximation
 property relative truncation and regularization (\cite{Schwartz57}). Thus, analogously to a normal operator, we can choose a compact (regularizing) completion $H \pm X$ with X compact and  
 $I(H \pm X)=H \pm X$.

Consider according to Oka (\cite{Oka60}), given $\varphi$ real and continuous, a characteristic surface through a point P 
according to $\sigma=\{ \varphi(x) > \varphi(P) \}$, that can be given according to
$\sigma=\{ x_{n}=Q(x_{1},\ldots,x_{n-1}) \}$, where Q is polynomial. The construction gives a
transversal intersection. Characteristic for strict pseudo convexity, is that
we have existence of a hyper surface $\sigma$ in nbhd P, such that $\sigma \subset nbhd P \backslash P$.

Weak solutions to the heat equation $(\delta_{t} + L)u(t,x)=0$, are given by $u(t,x)=e^{-t \Delta} f$, 
with $u \rightarrow f$, as $t \rightarrow 0$. Gårding (\cite{Garding60}) extends analyticity for weak solutions,
to vector valued Banach-spaces.

Analytic vectors: $b \in B$ is analytic, if $\mathcal{R}(x)b$ is analytic in x.  
Starting from $\mathcal{R}(x)b$ analytic, we consider $<b,d \mathcal{R}(x)> \in \mathcal{H}$. 
Consider $\mathcal{R}(U^{-1} V U)=\mathcal{R}(V)$,
where $\mathcal{R}$ is a representation invariant for equivalences. Note that $\mathcal{R}(U)b=\int b d U=<b,d U>=0$ in $\mathcal{H}$, 
does not imply $\mathcal{R}(U)b \rightarrow U$
continuously.  Assume that $(U,V)$ are conjugated and that $\mathcal{R}(U,V)$ constitutes a disk (positive measure). 
Then, we have given $f$ is harmonic
over $\mathcal{R}(U,V)$, that the (arithmetic) mean $\mathcal{R}(U,V)(f)$ is constant. Assume that $\overline{\mathcal{R}}(U,V)=\mathcal{R}(-V,U)$ 
and that $f \bot d \mathcal{R}$ iff $f \bot - d \overline{\mathcal{R}}$. Let $f \sim \frac{d F}{d \overline{\mathcal{R}}}$, 
then we have that $f \bot d \mathcal{R}$ iff 
$F$ is analytic in $\mathcal{R}$, that is $F \rightarrow \mathcal{R}$ is continuous. Note that $L^{1}$ can be approximated
by H over a strictly pseudo convex domain. 
Note also that $\mathcal{R}$ analytic with finite Dirichlet integral and single valued over a Riemann surface (U,V), 
implies  $\mathcal{R}$ linear over $(U,V)$. 

\subsection{Local lifting principle}
Assume that $\Omega$ is a domain of holomorphy, according to $d U=0$, with $bd \Omega$ such that 
$U=I$. 

\newtheorem{lyft}[prel]{Lemma}
\begin{lyft}
Consider $(x,y,z) \rightarrow (u_{1},u_{2})$
a planar domain. Existence of $F(d U_{1})=d U$, that is a lifting principle according to $d U=\rho d U_{1}$,
where $\rho$ is analytic, gives that $\Delta U = \rho \Delta U_{1}$. Assume that given $1/\rho \rightarrow 0$ in 
$\infty$ regularly,
$d V=\rho d U_{1} + \frac{1}{\rho} d U_{2}$ and
$d U_{2} = \vartheta d U_{1}$, we then have  that dV can be chosen as reduced in $dU_{1}$, when $\vartheta \rightarrow 1$ 
regularly. 
\end{lyft}
More precisely, given $U=I$ defines $bd \Omega$, we have in $\mathcal{L}_{ac}$ that $\Omega=bd \Omega$.
In $\mathcal{L}_{c}$, we include the possibility of $\vartheta$ const.
Assume that H defines a planar domain $(U_{1},U_{2})$, that is a cylinder and trivially an algebraic polyhedron. 
This gives a lifting principle. Assume that $\mid \vartheta \mid \leq 1$ regular and $\vartheta \rightarrow 1$ in $\infty$, 
we then have $d V=0$ iff $d U_{1}=0$ on compact sets.

Example : Assume that F is a figure and $G_{0} \subset G$ movements, such that $G_{0}$ preserves F,
that is consider $\mathcal{H}(G_{0}) \subset \mathcal{H}(G)$. 
For instance translation on the cylinder web C,
preserves C, when it is parallel with the z-axes, otherwise not. 

When $\mathcal{H} \subset \mathcal{K}$
the restriction is sub-nuclear. We assume here, that the restriction preserves disk neighborhoods. 
Starting from $Uf(\phi)=\int f(\phi)d U$, we must assume 
that the movement has a nuclear representation. When $\phi \in \mathcal{D}_{L^{1}}$, with relatively
compact sublevelsets, 
we can consider $f * \phi \in L^{1} \cap C^{\infty}$, that is $f \in \mathcal{D}_{L^{1}}'$. 
In particular, when $d U = \rho d U_{1}$, $\mid \rho \mid \leq 1$ defines a neighborhood of 0.
Nuclear spaces have the approximation property . Given $d U=\sigma d U^{\bot}$, where $d U^{\bot}$ 
is constructed relative 
a scalar product and is maximal with respect to convexity (\cite{Riesz37}), we then have where 
$\sigma = const$, an irregular approximation property. In this case, $R(d U^{\bot})$ (range) is not strictly 
included in $R(d U)$. Note, that given existence of 1-dimensional movement $d N$, with $d U / d N > 0$ 
and simultaneously $d U^{\bot} / d N < 0$, we see that $\sigma=const$ must be trivial.

Contact transforms preserve order of zero's, but given a locally convex topology, we consider the completion
of a discrete set to a compact convex set. 
Example: assume that $\Gamma$ a discrete set in the finite space, and $\tilde{\Gamma} = Cl(\Gamma)$ 
compact. Given d is the distance function to $\Gamma$, then d is locally 1-1, the same thing does not hold 
for $\tilde{d}$, that is given $d \mu$ locally 1-1 on $\Omega$, we do not necessarily have that $d \mu$ 
is locally 1-1 on $Cl(\Omega)$.

\subsection{Uniformity}

The order of a neighborhood is determined by the number of movements $U_{i} \in G$,
that generates the neighborhood, that is we assume locally in $\mathcal{H}$ that $\Sigma U_{i}=I$.  

Consider (U + V)f=f(u,v), where U,V are linearly independent. Assume that $(u,v)$ is dense in the domain 
for f, we can then write $f(u,v) \simeq f_{1}(u-v)$. Note that projectivity in the phase implies 
invertibility in the symbol space. 
Consider $(U(a),V(b))$, with $a + b=1$. Assume that the action is in the phase, U(a) V(b) f= $e^{(U(a) + V(b)) \phi}$,
that is over compact sets in a,b, we have that $U(a)=V(b)^{-1}$.
Density is dependent of the topology. Example: $\mathcal{H}_{2} \subset \mathcal{H}_{1}$, with 
$(U,U^{\bot})$ dense in $\mathcal{H}_{2}$,
assume $(U,U^{\bot},V^{j})$ is dense in the domain for intermediary spaces, that is ker $V^{1}$  
gives the polar to $\mathcal{H}_{1}$ relative $\mathcal{H}_{2}$. 
 
Assume that $I_{+}f=I_{-} f$, for instance $f \in C^{\infty}$ in $\mathcal{L}_{ac}$. The continuation
to $\mathcal{L}_{c}$ can be represented through symmetric neighborhoods, such that we have two-sided 
limits. Assume that the coefficients ($C^{\infty}$) can be approximated by polynomials $\mid P_{1} \mid \leq \mid \rho \mid \leq \mid P_{2} \mid$,
then a two-sided limit is taken as $\mid \rho / P_{1} \mid \rightarrow 1+$ and $\mid \rho / P_{2} \mid \rightarrow 1-$.
 
 Note that a necessary condition for $\big[d U,d V \big]$ to be absolute continuous, is that
 the movements applied to the coefficients, preserve absolute continuity. Note that given a 
 bounded domain for absolute continuity and monotonous movements, the movements can not be 
 iterated infinitely. Example : assume that $U+ (I - U)=I + V$,
 that is U is projective, where V=0. Determine W, such that $R(W)=N(V)$, we then have VW g=0 $\forall g$, 
 that is $V \bot W$, thus U is projective on R(W). Assume that $U$ is projective in the phase space, 
 we then have in the symbol space,
 $U(I-U)=I$, thus we have existence of $U \sim (I-U)^{-1}$. Sufficient for convergence, when 
 $U^{j}$ is absolute continuous for 
 $j \geq N_{0}$, is that $V(\Sigma^{N_{0}} U^{j})=0$.
 
 Consider symmetric neighborhoods: $(d U,d U^{*}) \rightarrow (d U^{*},-dU)$. Given symmetry,
 the corresponding differential form is pure, thus it has in the plane an analytic representation. 
 Consider $d U \rightarrow d U^{*}$ through conjugation, 
 it is then sufficient that $d U^{* *} = - d U$. Given $d U^{*} = \rho d U$, so that closedness is preserved,
 then we have over analytic domains, (deformations to) a harmonic representation.
 
 Assume that $\big[ X,X \big]=dV =0$, with $V=I$, this does not imply $\int X =I$. 
 In particular, given $d V=\rho X$, in the case where X is not one-parameter.
 By continuation to $\dot{B}$, we can define $d V$ in $\mathcal{D}_{L^{1}}'$, but $\dot{B}$ 
 is not dense in $\mathcal{D}_{L^{1}}'$. Note $X^{2}$ is absolute continuous and in $\mathcal{L}_{c}$ 
 which implies $X \in \mathcal{L}_{c}$, but does not imply absolute continuous.

\subsection{Generalized Sobolev spaces}

For a rigorous presentation of vector valued $L^{p}-$ spaces, we refer to (\cite{Schwartz57}).
Assume that $M^{o}=\{ d \mu \quad <f,d \mu>=0 \}$, linear measures to M, 
${}^{o}M^{o}=\{ f \bot d \mu \}$, that is given M a closed (Banach-) space we can write 
$M=R(d \mu)^{\bot}$, with $d \mu$ linear.

Consider $L^{p}(E)$ defined by $U T* \varphi = (UT) * \varphi \in L^{p}$, given the translates are dense 
in the domain. Assume that $d U = \rho d I$, where $\mid \rho-1 \mid \leq 1$ defines a closure, in particular 
when we assume $\rho$ linear over $(U_{1},U_{2})$. Note that given the topology for Banach-spaces and $\rho$ monotonous and 
analytic, with finite Dirichlet integral over a dense planar domain, we can give $\rho$ as linear. 
When characteristic sets for $L^{p}(E)$, are given by $U-I$, they can be considered in $L^{1}(E)$, for instance
$g^{p}=\rho f$, with $f \in L^{1}$. Assume that $T(\Omega) = \cup_{j} L_{j}$, 
where $L_{j}$ are radii. Assume on $L_{j}$, that 
$\mid \vartheta_{1} \mid \leq \mid \rho \mid \leq \mid \vartheta_{2} \mid \leq C$, 
where $\vartheta_{1},\vartheta_{2}$ are algebraic and corresponding to the same type of movement. 
This corresponds to
 $T(\Omega_{\vartheta_{2}}) \subset T(\Omega_{\rho}) \subset T(\Omega_{\vartheta_{1}})$.
The type of movement can be determined by taking $\underline{\lim}_{\rho \leq \vartheta_{2}}$ and 
$\overline{\lim}_{\vartheta_{1} \leq \rho}$, when the limits are equal.

Assume that $\Omega$ a domain for $d U=d V=0$ simultaneously. Assume that 
$P \in \Omega$ and consider a line L through P. If P is not isolated on L,
that is for instance dU of bounded variation =0 on L, we have that d U=0 can be continued to a disk (convex). 
Thus we can identify the movements in a pluri-complex sense. Consider $d U(x,y,z) \rightarrow d U_{1}(1/x,1/y,1/z)$.
Assume that $U \in G$ and $U_{1}$ defined on $R(G)^{\bot}$ $L^{p}-$ compact with an approximation property.
This means that $d U_{1}$ can be chosen as algebraic in $0$. Consider $d U_{1} = \rho d U$,
in particular when $\rho$ is linear, we can chose $d U_{1}$ as  analytic in the $\infty$. 

\subsection{P-convexity}

Assume $d U \in \mathcal{G}$, where $\mathcal{G}$ has maximal order, for the symbol f.
Assume that $<p d U(f),d N>=0$ with d N one-parameter and ${}^{t} p d N \neq 0$, for a polynomial p, 
that is given $d U(f)=0$ defines a strictly pseudo convex domain, then ${}^{t} p d N$ defines a normal. 
In particular assume that $\xi'=(\xi_{2},\ldots, \xi_{n})$ and $\xi_{1}=\xi_{1}(\xi')$ with algebraic dependence.
Then we have $\xi' \in \Omega$ implies $\xi_{1}=0$, that is $\xi_{1} \bot \Omega$.

We can compare with the construction of principal operators in Riemann geometry (\cite{AhlforsSario60}).
Assume that f is a symbol for a differential operator, so that we have existence of d N solvable, 
with $Nf \neq 0$ and $d U \bot {}^{t} p d N$ on a boundary $\beta$, i analogy
with the second principal operator. 
Assume further that, over a boundary $\alpha$, we have vanishing flux, that is $\int d U^{\diamondsuit}=0$.
Assume further that over $\alpha$, we have a 
``projective'' decomposition, so that $(U,U^{\bot})$, generates the domain, with $d U^{\bot}=\rho d U^{\diamondsuit}$, 
with $\rho \neq 0$ on $\alpha$
and where thus $U = I$ iff $U^{\bot}=0$ over $\alpha$. A necessary condition for a hypoelliptic symbol f, 
is that $U^{\bot}f=0$ implies $\alpha$ in the finite space. The two principal conditions, with this  
addition, defines in this article P-convexity (cf. \cite{Duflo}). 

\newtheorem{Pkonvex}[prel]{Proposition}
\begin{Pkonvex}
A P-convex group G, where $R(G)^{\bot}$ $L^{p}-$ compact 
with the approximation property, preserves hypoellipticity.  
\end{Pkonvex}

Assume that $U=I$ over $\alpha$ iff ${}^{t} U =0$ over $\alpha$. Assume that $\beta=\beta(\alpha)$ according to
 $d U = \rho d N$, with $\rho=0$ over $\beta$ and N normal. Further $d {}^{t} U=\rho' d {}^{t} N$ 
 with $\rho'=0$ over $\beta'$. 
 Given $(d U)^{\bot} \simeq (d N)$, we can write $N \bot R(dU)$,
 that is where R(dU) dense, we must have $N \sim 0$. Given $d N=0$ iff $d {}^{t}N=0$, we must have 
 $\alpha \cap ker N=\emptyset$ and $\alpha \subset \beta'$, that is on $\alpha$, U is only dependent 
 of N.  

Note that through the condition on 
$L^{p}-$ compactness for N, $<d U,d N>(f)=<d U,N(f)>$ defines a regular continuation. 

We note that $L^{1}(\Omega) \subset L^{1}(\Omega')$, given $\Omega' \subset \Omega$, that is given
$\Omega$ is a maximal domain for regularity (corresponding to a minimal set of invariants), the restrictions
are regular.  The boundary to a maximal domain for regularity is assumed to be given by a Riemann surface 
$\mathcal{W}$, that is $d U=\rho d I$ with $\rho \in \mathcal{H}$ implies $\rho = const$ on $\mathcal{W}$.

More precisely, $d W_{T}= \vartheta_{T} d W$ with $\mid \vartheta_{T} \mid \leq 1$, we then have 
$\mid d W_{T} \mid \leq \mid d W \mid$. Given $d W= \sigma d U_{1}$ with $d U_{1}$ harmonic, we have
further $\mid d W \mid \leq \mid d U_{1} \mid$.
The condition $<\rho d U,d W>=0$, when ${}^{t} \rho d W \neq 0$ implies a non-trivial kernel to d U,
that is a continuation (orthogonal) of dU.
Assume that $<\rho d U, d W>=0$ and assume that $d U \rightarrow d I$, we then have for $\phi \in L_{ac}^{1}$, 
that $<d I, {}^{t} \rho d W>(\phi)=
<\phi, {}^{t} \rho d W>=0$ and given $\rho$ polynomial so that $W(\rho \phi)=0$ $\forall \phi \in L_{ac}^{1}$, 
we have $W \sim 0$. Further, given $\varphi \in L^{1}$ with $\sigma \varphi \in L_{ac}^{1}$, given 
$\sigma \rho$ polynomial, the result can be continued.
Continuation through $p d I$, where p polynomial, defines locally a domain of holomorphy.
Assume that $\Gamma=\{ W-I \}$. Given W-I has an algebraic symbol, 
(for instance $R(W-I)$ has the approximation property) 
then $d_{\Gamma} \mid W-I \mid < C$, where $d_{\Gamma}$ denotes the distance to the boundary, 
locally defines a domain of holomorphy. 

Assume that $d U=\rho d U_{1}$ and $d V=\vartheta d U_{2}$, so that $ d W=\rho d U_{1} + \vartheta d U_{2}=0$,
with $\rho,\vartheta$ constants, implies $\rho=\vartheta =0$. Further, where $\rho,\vartheta$ constants, $\Delta W=\rho \Delta U_{1} + \vartheta \Delta U_{2}$, 
that is W is harmonically independent of $U_{1},U_{2}$, if $\Delta W$ is independent of
$\Delta U_{1},\Delta U_{2}$. Conversely, given W is harmonic and $\overline{\rho},\overline{\vartheta}$ 
absolute continuous, so that $\overline{d} \rho=\overline{d} \vartheta=0$ 
implies $\rho,\vartheta$ const, then when W is harmonically independent of $U_{1},U_{2}$, 
$\Delta W$ is independent of $\Delta U_{1},\Delta U_{2}$.
In the same manner, given $\rho,\vartheta \in L_{ac}^{1}$, we have that W is linearly independent of $U_{1},U_{2}$
iff W is symmetrically independent of $U_{1},U_{2}$.

 \section{A two-mirror model}

\subsection{Presence of a local pseudo base}
 
Assume that $u \in L_{ac}^{1}$ and $d \mu=u d I$. We then have that $d u=0$ implies u=const.
In $L_{c}^{1}$, there are examples of u such that d u=0 does not imply $u=const$, 
in this case $\Omega(u d I) \neq \Omega(d I)$ (\cite{Riesz56}).
$d I \in \mathcal{G}$ implies precense of a local base in $\mathcal{G}$.
In analogy with Oka's property: assume that $d U_{ij}=d U_{j} / d U_{i}$ on $\Omega_{i} \cap \Omega_{j} \neq \emptyset$
Assume that $0 \neq \rho_{ij} \in C^{0}(\Omega_{i} \cap \Omega_{j})$ and $\rho_{ij}=\rho_{ji}^{-1}$ 
further $\rho_{ij} \rho_{jk} \rho_{ki}=1$ on $\Omega_{i} \cap \Omega_{j} \cap \Omega_{k}$ and 
$\rho_{ij}d I=d U_{ij}$. Thus, the topology $\mathcal{H}$ must be such that $\mathcal{H} \ni \rho \rightarrow 1/\rho \in \mathcal{H}$
or $\rho=const$. Note that, given $\rho$ is absolute continuous (here in one parameter), 
where $\rho'' \neq 0$, a sufficient condition for $1/\rho$ to be absolute continuous, is that 
$ d \log \rho' \leq 2 d \log \rho$.
 
I analogy with Runge's property, if $d U=\alpha d U_{1}$ with $\alpha \rightarrow 1$ regularly, 
assume that $\alpha' , \alpha ''$ are sub sequences that both tend to 1 uniformly, we then do not necessarily have 
that $\alpha' \sim \alpha''$.
Assume for this reason $d U=\alpha d I$, where $\{ \alpha_{k} \}^{h} \sim \{ C_{k} \}^{h}$, 
and where we assume $\alpha_{k}$ entire. Note that $\{ \alpha_{k}=C_{k} \} \rightarrow \{ \alpha_{k}= C_{1} \}$, 
when $C_{k} \rightarrow C_{1}$.
Thus, given for instance a regular covering, we must have $\alpha_{k} \rightarrow \alpha_{1}$ regularly. 
In general, $X(\beta v)=v X(\beta) + \beta X(v) =0$ and $X(v)=0$ 
does not imply v=const.

 Starting from $d U = \rho d I$, we can define $\{ \rho=\lambda \}$ with leaves $(I)_{\lambda}$.
 Tangents are given according to $d U=\lambda d I$. 
 Note that $Uf=\int f d U=\lambda \int f d I=\lambda I f$,
 that is f can be seen as an eigenvector. 
 When (d U, $d U^{\bot}$) is dense in $L_{ac}^{1}$, then $d U + d U^{\bot}(\phi)=0$
 implies that $\phi=const$. 
 
 Starting from $d U \in \mathcal{G}_{q}$, there are n-q independent solutions. Assume that 
 $d U$ in a linearly independent base of analytic movements in $\mathcal{G}_{q}$, then given existence of $d V \in \mathcal{G}_{n-q}$, 
 with $dV = d U$ in the $\infty$, then $d V$ must be not analytic. Example: $d U=\rho d I$ with $\rho=const \neq 0$ on W, with W compact, but $d U \neq d I$ in the $\infty$, 
 that is U has no trace in the $\infty$.
 
 Consider $\int f dU=\int f(u)d \varphi(u)$, where $d U$ one-parameter.
 According to Riesz (\cite{Riesz23}) assume that $\varphi$ non-decreasing and
 u in a bounded domain. Consider $u(v)$ and $f(u)=h(v)$ and $\int f(u) d \varphi(u)=\int h(v)d v$, 
 where $h(v)$ is indetermined for $v=const$, these segments constitutes a countable set, that does not 
 influence the value of the integral. Note that we assume $d V$ one-parameter. This defines an orthogonal
 movement on compact sets.

  \subsection{Topological type for the movement.} 
 
 The topological type for movement $d U = \rho d I$, is given by the topology for $\rho=const$ on 
 $\Omega$, that is the least fine topology, so that $\Omega$ is invariant. If the invariant movement can 
 be continued by continuity, then $\Omega$ is not minimal. More precisely, starting from $d U \in E'$,
 where $E'$ has the approximation property (by truncation) and
 $<T,d U>(\phi)=<T(\phi),d U>$, we must have $E^{'*} \ni T(\phi) \rightarrow \rho^{*} T(\phi) \in E^{' *}$. 
 When $E^{'*}$ has the approximation property (by truncation), we must have $\rho^{*} \rightarrow 1$ iff $\rho \rightarrow 1$.
 Example: consider $U \rightarrow {}^{t} U$ as an extension, that is 
 $D (U) \subset D(U^{t})$, we then have ${}^{t} \Omega \subset \Omega$. However, the movements
 do not have to be of the same type.
 
 Example: assume that $d y / d x=- \eta / \xi$, that is $y=y(x)$, given $\xi$ divides $\eta$. 
 Given $d y/d x=0$, we must not have that y is constant anywhere. Given 
 $d U_{1}$ is absolute continuous, but $d U_{2}$ continuous,
 we have that $\Sigma_{1} \simeq \Sigma_{2} = \{ d U_{2}=0 \}$ does not imply $\{ U_{1}-I \} \simeq \{ U_{2} - I \}$
 
 \newtheorem{sng}[prel]{Lemma}
 \begin{sng}
 Assume that $d U = \rho d V$ is a continuation of dV and consider $T(\phi) \in C^{0}$ 
 (for instance regularizations). Thus, $<T(\phi),\rho d V>=<\sigma T(\phi), d  V>$, where we assume 
 ${}^{t} \sigma : \dot{B} \rightarrow \dot{B}$. We see that $d U$ continues invariant sets, 
 where $\sigma=const$. Assume that $\{ {}^{t} \sigma=const, d {}^{t} \sigma=0 \}$
 algebraic, for instance $\log {}^{t} \sigma \in L^{1}$, then there are no non-trivial continuations 
 of $dV \neq 0$, on this form, that preserves symmetric forms. In the same manner, given 
 $\{ {}^{t} \sigma=const, \overline{d} \sigma=0 \}$ algebraic,
 there are no non-trivial continuations of $dV \neq 0$, on this form, that preserve harmonicity.
 \end{sng}
  \subsection{$L^{p}-$ compact complement}
 
 Given P,Q are combined through a transformation in $\mathcal{G}_{r}$, simply transitive (reciprocal),
 we can completely define $\mathcal{G}_{r}$, by letting P,Q vary. 
 Assume that P,Q are combined by an arc $\gamma \in C^{1}$ and $I \in G$,  given $d U_{\rho}(f) \equiv 0$,
 when $\rho \rightarrow 1$, where $d U_{\rho}=\rho d U_{1}$, then $G$ is completely defined by $\rho$,
 when $P,Q$ varies.
 Assume that $\gamma(P)(t)=Q(t)$ and form trajectories between P,Q, that are generated by $\mathcal{G}_{r}$, 
 that is $R(\mathcal{G}_{r})$. We then have that $R(\mathcal{G}_{r})$ is locally convex relative 
 $\mathcal{G}_{r}$. A sufficient condition for $\mathcal{G}_{r} \simeq \mathcal{H}_{r}$, is that 
 $\Omega(\mathcal{G}_{r}) \simeq \Omega(\mathcal{H}_{r})$ (a projective mapping).
 
 $R(G_{q})^{\bot}$ $L^{1}$-compact with the approximation property, corresponds to the proposition
 of a regular, compact complement.
  Given f analytic, we have that $\{ d \mid f \mid \leq c \}$ is a domain of holomorphy,
 further $m(\{ d \geq R \}) \rightarrow 0$, as $R \rightarrow \infty$, that is Lebesgue compact means 
 $\{ \mid f \mid \leq 1/d \leq 1/R \} \subset \subset E_{n}$
 
 Hypoellipticity has a pluricomplex definition, where a necessary condition is that no lines, contribute 
 microlocally, that is we do not have non-trivial lineality.  We can for this reason consider $\mathcal{G}_{1}$, corresponding to a line through 0, 
 for instance a reflection axes. 
 Thus, to discuss the necessary condition for hypoellipticity in $L^{1}$, it is sufficient
 to consider reflection (involution). 
 
 Given E is separated and $F \subset E$ with $I \in \mathcal{L}_{c}(E,F)$ (restriction) and 
 where F has the approximation property, we have that E has the approximation property  (\cite{Schwartz57}). 
 Assume for this reason that $E'$ is given by dV $\in \mathcal{G}$ and
 $d U=\rho d V$. In weighted $L^{p}-$ spaces, a necessary condition for the inclusion $F \subset E$
 is that $1/\rho^{*} \rightarrow 0$ in $\infty$.
 In particular, given existence of $\beta$ with $\rho=0$ implies a strict inclusion $F \subset E$.
 Further, when $d U \rightarrow d I$ regularly, we have $d V \rightarrow  d I$ regularly.
  
\subsection{Dependence of orientation}

For a projective representation, we have that $U(I-U) \simeq 0$, given uniformity, $R(I-U) \subset \mbox{ ker} U$ or
$R(I-U) \cap R(U) = \{ 0 \}$.  
Assume that $d U =\Sigma \rho_{j} d U_{j}$, where $d U_{j}$ one-parameter. Where $\rho_{j}=const$, we have that d U are linearly dependent 
of base vectors, where $d U_{j}$ are analytic.
Further, $d^{2} U = \Sigma d \rho_{j} d U_{j} + \rho_{j} d^{2} U_{j}$, thus given $d U_{j}$ 
convex (absolute continuous), we have where $\rho_{j}=const> 0$, that U is convex as a movement.

\newtheorem{konvex}[prel]{Lemma}
\begin{konvex}
For $d U = \rho d V$, we have existence of $W_{2} \subset R(d V)$, so that d U is convex over $W_{2}$. 
\end{konvex}

More precisely, sufficient given $d U=\rho d V$, is that $\frac{d \rho}{d v} \geq 0$ over $W_{2}$.
Consider $(x,y,z) \rightarrow (v) \rightarrow \rho(v)$. Given $\rho \in C^{\infty}$, we can 
obviously determine $W_{2}$ with $\frac{d \rho}{d v} \geq 0$ on $W_{2}$. Given $Uf=g$, that is $\int \rho(v) f(v) d I(v)$,
and consequently $Uf = V(\rho f)=g$ or $\rho f(v)=g(v)$, why we must assume $g \in R(V)$.

\newtheorem{harm}[prel]{Lemma}
\begin{harm}
Assume that the coefficients to X are given by f, $X=X_{f}$, given $\{ f,g \}=0$ that is 
the coefficients can be given as results of reflection,
then we have existence of v harmonic, with $\big[ X_{f},X_{g} \big](v)=0$.
\end{harm}
More precisely, $\big[ d U_{1},d U_{2} \big]= (d U_{1} (\xi^{2}) - d U_{2}(\xi^{1})) v_{x} + 
(d U_{1}(\eta^{2}) - d U_{2} (\eta^{1})) v_{y} + (\xi^{2} \eta^{1} - \xi^{1} \eta^{2}) v_{xy} + (\xi^{1} \eta^{2} - \xi^{2} \eta^{1}) v_{yx}$,
that is we have existence of v symmetric, so that the bracket is 0, thus we have existence of v 
harmonic. The condition $\{ f,g \}=0$ makes the two first terms disappear.

Jacobi can be seen as a two-mirror reflection model: $0 = - \{v, \{ f,g \} \}= \{ f, \{ g,v \} \} + \{ g, \{ v,f \} \}$.
For instance, when f,g are reflection results of v. Given existence of v, so that $\{ f, \{ g,v \} \}= \{ g, \{ f,v \} \}$,
that is we consider $v \rightarrow g \rightarrow f$ as equivalent with $v \rightarrow f \rightarrow g$,
we have a double base model.

Assume that U is reflexive, but not projective according to $U^{\bot}=(I-U) + V$ and $R(W)$, so that 
$V=V^{\bot}$. Consider a regular representation $K$ on $R(W)^{\bot}$, corresponding to $N(W)$ 
(through the projection mapping), 
when we assume $W_{T} \rightarrow 0$, we get a corresponding very regular representation 
$K_{T}(V,V^{\bot})$.

\subsection{Envelop representation}

Define $L^{1}(dU)=\{ g / \rho^{*} \quad g \in L^{1}(d I) \}$. Assume further that $\mathcal{G}$ can be 
represented through $\mathcal{M'}=
\{ \rho \quad \rho d I \in \mathcal{G} \}$. We can then construct $L^{1}_{c}(d I) \simeq \rho^{*} L^{1}_{ac}(d I)$ 
Note that $g=\rho^{*} f \in L^{1}_{ac}$
does not imply $f \in L_{ac}^{1}$.

More precisely, for $\mathcal{L}_{ac}$, Id has support in the entire space and $\mid \rho -1 \mid \leq 1$ 
defines a convex neighborhood in $\mathcal{M'}$, that is given $f \in L_{ac}^{1}$, we assume $\rho^{*} f \in L_{c}^{1}$ 
and we can generate a continuation.

Where $f \in L_{ac}^{1}(d U_{j})$, the invariant sets $\Omega_{j}$, are zero-sets to $d U_{j}$ and through uniformity,
for the continuation to $L_{c}^{1}$, we have $\Omega_{j} \subset ker (d U_{j})$ or 
$\Omega_{j} \cap ker (d U_{j}) = \{ 0 \}$. The continuation to $L_{c}^{1}$
is not reducible, for instance given $\Omega=\bigoplus \Omega_{j}$ implies 
$d U = \Sigma d U_{j}$ and when $d U \neq 0$, consider $d U_{j} = \rho_{j} d U$ with 
$\rho_{j}=1$ on for instance $\Omega_{j} \cap \Omega_{k}$, $j \neq k$.

For a $\rho-$ stable subspace $W \subset V$ (\cite{H/S}), we have that $W \cap (\Sigma_{D \in \Omega} V_{D}) = \Sigma_{D \in \Omega} W \cap V_{D}$.
For instance, given $\Sigma \mid I_{j} \mid < \delta$ implies $ \Sigma \mid \int_{I_{j}} f d U \mid < \epsilon$,
that is $\int_{I} f d U=\Sigma \int_{I_{j}} f d U$, for f fixed. In $\mathcal{L}_{ac}$ 
semi-simple algebras are $\rho-$ stable.

Assume $\pi : \mathcal{G} \times B \rightarrow B$ continous, where B a Banachspace with I. 
Assume that $d U=\rho d I$. The problem $T(\pi W) \subset \pi T(W)$ is for this reason dependent of
 $\rho \rightarrow {}^{t} \rho$, that is $\rho d I(f)=0$ implies $d I ({}^{t} \rho f)=0$. 
The problem is obviously solvable given $\rho$ algebraic, in particular linear. 

Note X a nuclear Banach space, implies $X'$ a Banach space.
Consider $(AC)'$: given f absolute continuous, we have $\int f d I=\int f' d x=f$.
Assume that $d U=\rho d I$, with $\rho f$ absolute continuous, it is sufficient that $\rho$ is absolute continuous. 
In this manner, we identify $I({}^{t} Uf) \simeq \int \rho f d I \simeq UI(f)$.

According to Nullstellensatz, $\mid AB \mid \leq c$ implies $\mid B \mid \leq c'$, given A reduced.
However AB absolute continuous, does not imply B absolute continuous.
Consider $L : R(AB) \rightarrow R(B)$ sub-nuclear, given $I$ sub-nuclear, we have $R(B)$ nuclear. 
Given A reduced and not $I$, for instance $I=\delta_{x}$ in $\mathcal{L}_{c}$, we have $AB=I$ implies $B \neq I$, 
that is we have that $R(AB-I) \rightarrow R(B-I)$ is not sub nuclear.
When we consider the same problem modulo $C^{\infty}$, we have that ker $B \subset C^{\infty}$, 
that is $AB-I \in C^{\infty}$ implies $B-I \in C^{\infty}$.

\subsection{Invariant functionals}
For $T \in \mathcal{H}(E)$, we can define a support for T, with respect to $E'$ (\cite{Schwartz57}). 
Assume that we have existence of $d U_{j}$ with $U_{j} T \in \mathcal{E}'$, then we do not necessarily 
have $T \in \mathcal{E}'$. However, if we have existence of ${}^{t} U$ locally 1-1 (modulo $C^{\infty}$), 
that is ${}^{t} U \phi \in C^{\infty}$ implies $\phi \in C^{\infty}$, we then have $T \in \mathcal{E}'$.
Given E is a Banach space, the difference is not significant. Note also that $d U$ has support in its domain, 
that does not have to be compact.

Assume that $<f,d U>=0$ in $\mathcal{H}_{c}$, with an approximations property, implies
$f \in \mathcal{H}$, that is scalarly regular is interpreted as with regular orthogonal.
Given $\rho f \in \mathcal{H}$, as $\rho \rightarrow 1$, we must have $f \in \mathcal{H}$. 
When $\rho f \in \mathcal{H}$ implies $(1/\rho) f \notin \mathcal{H}$, then $\rho$ must be constant.
Further, assume that $T \bot d U$ in $(L_{ac}^{1})'$, where T is defined on a pseudo convex domain,
that is with a boundary of order 0, assume planar. Given $f \in L_{ac}^{1}$,
we have $Uf=f$ over the support to T, since I is 1-1 over the boundary and we have presence of a normal model.

The problem to determine a restriction $\mathcal{L}_{c} \rightarrow \mathcal{L}_{ac}$ can be compared with 
presence of a separating functional, For instance the radius of convergence. Assume that $R(d U,d U^{\bot})$ dense in 
$\mathcal{H}_{ac}$. When we consider $\int f d U \rightarrow \int f d I$, it is necessary
that $d U$ is monotonous, that is we consider d I one-sided. Consider a complexification, where $(d U, d U^{\bot}) \rightarrow (d I_{+},0)$
iff $(d U, d U^{\bot}) \rightarrow (0,d I_{-})$, then Id can be continued uniquely, to a separated, quasi-complete
space $\mathcal{H}_{c}$ (\cite{Bourbaki64}).

Completing sets (\cite{Schwartz57}): Given E quasi-complete, we then have 
for every bounded set $B \subset E$, that B is completing.

Assume that L is a line through 0, consider $B' \cap L = \emptyset$ iff $\{ d V < \lambda \} \subset \subset E$,
and consider the continuation $B' \rightarrow B$. For instance when $B \cap L$ is a closed segment,
we have $\{ d V < \lambda \} \rightarrow \infty$, that is non-compact. Thus, given $d V < \lambda$
over B, there is a maximum-principle and clustersets.

Example: assume that L a radius as above and $L \cap \Omega(U)(d f) = \emptyset$, for a fixed symbol 
$f \neq 0$, but $\Omega(d U)(f) \cap L$ a closed segment, in particular when $U$ is translation, then 
 L contributes microlocally to f, but not to df.

Assume that $G \ni U \rightarrow U^{\bot} \in G$. 
Through a change of variables $(x,y,z) \rightarrow (u,u^{\bot})$, we see that
$\mathcal{H}(G^{*}) \rightarrow \mathcal{H}(G^{*})$ has sub nuclear identity $u \rightarrow u^{\bot}$ 
on a disk neighborhood, for symmetric distributions.

Regarding scalar invariance, assume that $e' = d I$ and $<T,d I> \simeq I T$, that is TI=IT implies $\int T(\varphi) d I = T(\varphi)$
in the sense of functions, given that $T(\varphi)$ is absolute continuous and in the sense of distributions, 
when $T(\varphi)$ is only continuous.  
Note $I \notin \dot{B}$ and consider $d U=\rho d I$,  when $\rho$ is algebraic in x,
then $\rho$ maps first surfaces on first surfaces, but given that $\{ \rho=const \}$ are non-trivial, U may only be 
defined in the sense of distributions.
Let $M_{\rho}=\{ f \in \dot{B} \quad \rho f \in L_{ac}^{1} \}$. Obviously it is also  $L^{1}$-compact. 
Given $<T,dU> \in \mathcal{D}_{L^{1}}'$, we have then over $M_{\rho}$, $<T,dI> \in (L^{1}_{ac})'$
and $<T / \rho,d U> \in \mathcal{D}_{L^{1}}'$.

\subsection{The two-mirror model}

Consider the two-mirror model through $I \in G$,$J \in G^{\bot}$, $I \neq J$, but $I \rightarrow J$
projective. $I \in G$ implies $\exists U^{-1} \in G$, but we do not assume $J \in G$.
Given $\psi : L \rightarrow M$ with $\psi=AB$ and $AB=BA$ with $B^{2}=I$, we have $\psi^{2}=A^{2}$
and given A projective, $A^{2}=A$.

More precisely, consider $M' \rightarrow N' \rightarrow L'$ according to $d U = \rho_{2} / \rho_{1} d V$ 
and $d U=\rho_{2} d W$ and $d V=\rho_{1} d W$. $N'$ can be suitably chosen, for instance absolute continuous. 
Then d W(f)=0 implies Wf=f, but we do not necessarily have that $d U=d V=0$, implies $U=V=I$. Note that 
a necessary condition for inclusion of domains for respective measure, is that $\rho_{2} / \rho_{1} \rightarrow 0$ in 
$\infty$. Consider I locally 1-1 and J surjective, with ${}^{t} J \sim I$. Assume that $L'$ $L^{p}$- compact 
with an approximation property. Assume that $p \leq \rho_{j} \leq q$, where p,q polynomials. Consider relative 
the measures, that we have one-sided limits $\lim \rho_{j}/q \rightarrow 1-$ and $\lim \rho_{j}/p \rightarrow 1+$, 
when the measures coincide, the corresponding moment problem is determined. 
This means that the system $(G,G^{\bot})$ according to the above has an approximate solution. 
Note in the model as above, given $d U \in \mathcal{G}$, the inverse can be given as an approximate inverse.

Given a projective decomposition of d I by $d U_{1},d U_{2}$ linearly independent in $\infty$ and analytic, we have that for 
corresponding invariant sets, that $\Omega_{1} \cap \Omega_{2}=\emptyset$. Given linear dependence 
we have through uniformity, $\Omega_{1} \simeq \Omega_{2}$. For the restriction to a line, consider the scaling 
variable s, given linear independence
we have independence of s, that is $d U_{1} = d U_{2}$ on L, implies s finite. Given $d U_{1}/d U_{2}=const$, 
when s large, we have linear dependence. 

Example: Assume that rad U $\rightarrow U^{\bot}$ projective, for instance  rad dU $\simeq dU + c$. Assume that 
$\frac{d^{m} \Phi}{d t^{m}}=0$ implies $\{ \Phi=1 \}$ a discrete set. Further, $\frac{d \Psi}{d t}=0$, 
that is $d U^{\bot}(\Psi)=0$. Thus the set $\{ \Phi=1 \}$, so that $d U^{m}(\Phi)=0$, has the same order as 
$\{ \Psi=1 \}$. Assume that $d U=\rho d I$, it is then necessary that $\rho^{m-1} \Phi=\Psi$ in 
$H$. Assume that when the condition holds in $H \cap L^{1}$,
we have when $\Psi=1$ that $\Phi=1/ \rho^{m-1}$ in H and further given $\Phi=1$, we have that 
$\rho^{m-1}=1$ gives multivalentness for movement $U^{\bot}$.

\subsection{Conjugation}
Denote $F^{\diamondsuit}(\xi,\eta)(\alpha)=F(-\eta,\xi)(\alpha)=F(\xi,\eta)(\overline{\alpha})$. 
Assume that the domain for F is given by $(U,U^{\bot})$.
Assume that $d U$ has coefficients $(\xi,\eta)$, we then have  ${}^{t}d \overline{U}$ corresponds to 
$(-\eta,\xi)$, that is $d U^{\diamondsuit}$. Thus, given $d {}^{t} \overline{U} = \rho d U$ with $\rho$ analytic 
and given $d U$ is closed, then d U is harmonic.

Assume that $\rho \frac{\delta f}{\delta u_{1}} + \vartheta \frac{\delta f}{\delta u_{2}}=0$. 
Assume that $d U_{2}=\sigma d U_{1}$, so that $u_{2} = \psi u_{1}$, where $\psi$ continuous. 
The condition $(u_{1},u_{2})$ (linearly) convex, implies $u_{2} / u_{1} \rightarrow 0$ in $\infty$. 
The condition $\frac{d U_{2}}{d U_{1}} \rightarrow 0$ implies $\vartheta/\rho \rightarrow 0$ in 
$\infty$. Assume that $(d U,d U^{\bot})=(\rho d I, \frac{1}{\rho} d I)$ projective, with $\rho + 1/\rho$ reduced, 
this implies absence of essential singularities in $\infty$.

Assume that $W$ is chosen so that $d U + d U^{\bot}=(\rho + \rho_{1}) d W$ 
with $\rho + \rho_{1}=1$ on $R(d W)$. Assume that U algebraic on R(W), given $I \in R(W)$ we have existence of 
$U^{-1}$ on R(W).
Given U projective, we can construct V so that $U + V=I$ with $U=0$ implies $V=I$ and $U=I$ implies $V=0$, 
when U is not projective, assume that we have existence of V so that
$d L_{t}(U,V) \in \mathcal{G}$ $\forall t$ and $\frac{d L_{t}}{d t} \neq 0$, that is a (linearly) convex continuation 
of U to maximal order with $\Omega_{U} \cap \Omega_{V}=\{ 0 \}$ or $\Omega(UV) \subset \Omega(V)$. 
Assume that $U+V$ is algebraic and has maximal order in the phase, we then have $UV f \simeq Uf$ 
or UV has trivial invariants. 

Assume that $d V_{1} \leq d U \le d V_{2}$, and $d V_{1}=\sigma d U$, $d V_{2}=\rho d U$.
Thus we have $d V_{2} = (\rho / \sigma) d V_{1}$ and we assume $\mid \rho / \sigma \mid \leq 1$.
($\rightarrow 0$ i $\infty$). Let $d I_{-}$ correspond to $\sigma \rightarrow 1$. A sufficient condition
to determine d U, is that $I_{-}(U) = I_{+}(U)$.
More precisely, $\forall d U$ of bounded variation with $\int g d V_{j}=0$ $j=1,2$, 
we have that $\int g d U=0$.

 \section{Continuation}
Given $f(x,y)$ not identically zero, we have that $f-Id \in C^{\infty}$ implies $f \in C^{\infty}$ in $L^{1}_{ac}$, 
the correspondant in $L^{1}_{c}$ gives $f - Id \in C^{\infty}$ implies f very regular. 
In particular, $(U-U^{\bot}) f \in C^{\infty}$ implies $f \in C^{\infty}$, 
given $U^{\bot}=I+U$, that is projectivity.

Assume that $T \in (L^{1}_{ac})'$, we then have for $T \sim Id$, $(\mbox{ supp T})^{c}=\{ 0 \}$, that is $\mbox{supp} (T- Id) = \{ 0 \}$.
Consider $T_{j}=Id + K_{j}$, where $K_{j}$ corresponds to regularizing action, analogous to a stratification.
Assume that $\Omega(T)=\{ f \quad T(f) - Id (f) = 0 \}$, we then have for $f \in L^{1}_{ac}$ that 
$\Omega(Id)^{c}=\{ 0 \}$. Assume $R(K_{j})$ defines a neighborhood of $\Omega(T)$. When Id 
locally 1-1 in $L^{1}_{ac}$, we have that ${}^{t} Id$ is locally surjective in $(L^{1}_{ac})'$,
where $(L^{1}_{ac})'$ is defined so that every measure can be given by $d \mu \rightarrow {}^{t} \rho d \mu$, with $\rho L^{1}_{ac} \subset L_{ac}^{1}$, 
analogous to Radon-Nikodym. Given the continuation includes $\rho=1$, we have that $(L^{1}_{ac})'$ 
is nuclear. When $R(K_{j})$ $L^{p}-$ compact, with an approximations property, we have that 
$\{f \quad \mid T(f)- Id(f) \mid \leq \lambda \}$ are relatively compact sets. In particular $K_{j}(\phi)=\phi(k_{j},{}^{t} k_{j})$, 
where $K_{j} \rightarrow {}^{t} K_{j}$ preserves regularizing action, analogous to uniformity. 

Concerning P-convexity, assume that $\alpha$ a boundary in x,y, so that $\Omega(T) \simeq \{ T-I=0 \}$ 
and that we have relatively compact sub level surfaces in the convex
closure of $\Omega(T)$, that includes the boundary $\alpha$. The condition $d T/d N=0$ over a boundary 
$\beta$, under the condition $R(T)^{\bot}$ is $L^{p}-$ compact, means that $K_{j}$ have support on one side of a hyper plane $\supset \beta$. 
Under the condition that we have an approximations property in $R(T)^{\bot}$, we can
assume that $\Omega(T)^{c}$ is locally removable. 
 
 Note that $< Uf, U(\phi)>= <Uf, <d U_{x}, \phi> >=<<f, d U>, <d U,\phi>>$.
 Assume that $d U = \sigma d I$ and $<d U(f),\phi>=<f,(\eta_{y} + \xi_{x}) d I(\phi)> + <f,d U(\phi)>$ (we assume $\xi,\eta$ algebraic),
 that is given $(\eta_{y} + \xi_{x}) > 0$ and $<d U(f),\phi>=<f,d U (\phi)>$ $\forall \phi$, then $d I(f)=0$, thus $d U(f)=0$.
 
 Given $\int f d U=(U+V)f$, where V defines the polar, we then have $U+V \in G_{r}$ in $\mathcal{L}_{ac}$,
 but $U + V \in G_{r+m}$ in $\mathcal{L}_{c}$, that is the order can be affected by $\mathcal{L}_{ac} \rightarrow \mathcal{L}_{c}$.
 
  Example: assume that $\Sigma=\{ U=I \}$ and $\Omega=\{ d U = d I \}$, further $d U=\rho d U_{1}$.
 We then have $d U_{1}=0$ implies $d U=0$, but we may have $\Sigma_{1} \neq \Sigma$. Given $d U_{1} = d I$, 
 we have $d U = \rho d I$, that is $\Omega_{1} \supset \Omega$. This means that the type for $d U$ is different from $d U_{1}$, 
 that is given $d U,d U_{1}$ have the same type, 
 then $\Omega_{1}$ is not minimal.

 \subsection{Relative continuation}

Relative inverses: Assume that $X=R(A) \bigoplus R(B)$ (two-mirror model), in particular given 
$B=I-A$ and $A(g-h)=f-h$, that is given $f=g$, Ag=f iff Ah=h and A + B=I. 
Assume that $\psi(h)=g$, we then have $A(\psi(h) - h)=f-h$, that is $A\psi(h) - f=Ah-h$ and given 
$\varphi(h)= f$ has continuous inverse,
there is a solution to the equation. Assume for instance that $\psi \circ \varphi^{-1}$ surjective on $f \in R(A)$ and that 
$\mathcal{H}$ is locally convex.

Note that $U(dU^{\bot})=\int d U^{\bot} d U$ and $U^{\bot}(d U)=\int d U d U^{\bot}$.

Given $d U^{\bot} / d U=\rho$, we have $U(\rho d U) + U^{\bot}(d U)=U(dU^{\bot}) + U^{\bot}(d U)$.
Consider $U(\rho d U) + U^{\bot}(d U)= \lambda d U$. Given f Hamiltonian with right hand sides $(X,Y)$ 
and g Hamiltonian with right hand sides $\xi,\eta$, then through Lie, $dU=X_{g}$, must be one-parameter. We can 
compare with $-X^{*} Y + Y^{*} X=R \neq 0$ (contact transform),
where the equation is degenerate for $R=0$. We consider movements related to one-parameter movements according to
$\big[ d U,X \big]=0$ or $\big[ \big[ d U,X \big], Y \big]=0$.
Given $R(d U)$ compact, regular on $R(W)$, for instance $d UW=\rho dW$, with $\rho$ analytic and $\rightarrow 0$ 
in $\infty$, we then have that $N(d {}^{t} U)$ (kernel) defines invariants relative $R(W)$,
analogous with Fredholm's alternative; that is, under an involution condition, so that $U + U^{\bot}=I$, 
we have that $\varphi' \in N({}^{t} U)$ iff $\varphi' \bot UW f$, for some $f$.

\subsection{Non-linear continuation}
Assume now that  $\mathcal{L}_{ac} \rightarrow \mathcal{L}_{c}$ corresponds to 
$I \rightarrow I + K$, with $ I \sim \delta_{x}$ i $\mathcal{L}_{c}$. When K is a regularizing operator, 
$I + K : \mathcal{D}' \rightarrow \mathcal{D}^{' F}$.
Example:  $\Sigma U_{i}=I$ continued to $U_{1} + C^{\infty} + \Sigma U_{i}$,
that is $R(\Sigma_{2}^{N} U_{i})$ generates $R(U_{1})^{\bot}$ and the polar. 

Consider $\psi : I_{K}(f)(x) \rightarrow f$
and $\psi' : I_{K}(f)(x) \rightarrow x$. Given K=Id and the topology $\mathcal{L}_{ac}$, we have 
$\psi$ continous and if further f is reduced, we have $\psi'$ absolute continous (outside a compact set).

We have that $L^{1}$ over a strictly pseudo convex domain, can be approximated by H (\cite{Garding60}).
For a pseudo convex domain, given the normal can be approximated two-sided by polynomials, the corresponding
moment problem is determined. Note that
d U monotonous, means that $d U_{+} \rightarrow d I_{+}$, in particular given the coefficients can 
be given according to $d U_{-} \simeq \overline{d U_{+}}$, we have $I \simeq \mbox{ Re } I$. Note 
that given $f \rightarrow I_{K}(f)$ locally 1-1 i $L_{ac}^{1}$, we have $I_{{}^{t} K}(f)$ locally surjective
in $(L_{ac}^{1})'$, for instance $d U = \rho d I$ with $\rho f$ absolute continuous and $\rho \neq 0$. 
Further given $I_{K} + R=I$, where $R$ $L^{1}$- compact, 
gives $K$ as very regular in $L_{c}$, through the condition on the normal.

Consider $<d U^{\bot},\phi> + v(x)$, a non-linear continuation, that is discontinuous.
Note that non-linear in phase implies non-algebraic in symbol space, that is absence of regular approximation property.
Note that given $d U + d U^{\bot}=d I$, we have that $\Omega(d U)=N(d U^{\bot})$.
Given $R(dU^{\bot}) \cap \Omega(dU)$ is non-trivial, we have multivalentness and vi choose a 
non-linear representation for $d U^{\bot}$, that is $<d U^{\bot}(\phi),T> + <v(x)(\phi),T>$,
where T is for instance in $\mathcal{D}_{L^{1}}'$ and $v(x) \phi \in \dot{B}$.
Given $d U^{\bot} + d U =d I + d V$, the 
polar for reflexive U, is characterized by $d V=d V^{\bot}$.
 
 \section{The approximation property }

Consider $d U \in E'$ and $\mathcal{H}(E)$. When U is analytic over f, 
we have that $Uf \rightarrow f(u)$ is continuous, a property that is not necessarily preserved,
when we continue $U$ to $L^{1}$.
Density for $(u_{1},u_{2})$ in E relative $\mathcal{H}$, when $UT \rightarrow T$ in $L^{1}$, 
corresponds to the proposition that $L^{1}(u_{1},u_{2})$ is dense in $\mathcal{H}(E)$. 

Starting from Lie, given $R(dU)$ has constant order, the order for $d U^{\bot}=d I$ must be constant. 
Assume that $(d U,d U^{\bot})$ is a complete system, $d V=\rho_{1}d U$ and $d V^{\bot} = \rho_{2} d U^{\bot}$.
Assume that $d U^{\bot}=\vartheta d U$, we then have $d V^{\bot}=\rho_{2} \vartheta d U$, that is $\rho_{1}^{\bot} \sim \rho_{2} \vartheta$
and the conjugation can be preserved. More precisely, assume that $d UW = \rho_{2} d W$ and $UW=V^{\bot}$ and $W U^{\bot}=V$, $d W U^{\bot}=\rho_{1} d U$
$=\rho_{2} \vartheta d U = \vartheta d U W = d U^{\bot} W$. Thus, $U^{\bot}=W^{-1} U W$. 
The approximation property can be continued, given $\vartheta$ regular,
that is $d U \bot d U^{\bot}$ implies $d V \bot d V^{\bot}$, as long as $\rho_{1} - \rho_{2} \vartheta \neq const$.

Note that for the approximation property to be continued $F \rightarrow E$, where $F^{\bot}$
$L^{p}$- compact, the inclusion $F \subset E$ must be well-defined, for instance a condition equivalent
with $L^{p}_{\vartheta} \subset L^{p}_{\varphi}$, is that $\varphi / \vartheta \rightarrow 0$ in $\infty$.

Consider $<I_{T}(\varphi),\phi>=<I_{T},\phi>(\varphi)$, where $\varphi \in \dot{B}$ and $\phi \in \mathcal{D}_{L^{1}}$
that is $I_{T}(\phi) \in \mathcal{D}_{L^{1}}'$ and $I_{T}(\varphi) \in (\mathcal{D}_{L^{1}})'$. Consider $d U = \rho d V$, 
that is $\rho \rightarrow 0$ i $\infty$, thus the corresponding ideals are inclusive $(I_{1}) \subset (I_{2})$.
This means that an approximation property for $(I_{2})$ can be derived from the corresponding property 
for $(I_{1})$. Thus, for instance $\rho \in \mathcal{D}_{L^{1}}$, implies given $d U \rightarrow d I$ 
regularly, that we can approximate $d V \rightarrow d I$ regularly, for instance through $\rho \rightarrow 1$. 

Example: Assume that in $L_{ac}^{1}$, $d U^{\bot}=\rho d U$ and $p_{1} d V_{1}=\rho d U = p_{2} d V_{2}$, where $p_{1},p_{2}$ are 
polynomials and $p_{1} \leq \rho \leq p_{2}$, given $\rho / p_{2} \rightarrow 1$ we have $\Omega(d U)=\Omega(d V_{2})$ 
and in the same manner, when $\rho / p_{1} \rightarrow 1$,
we have $\Omega(d U)=\Omega(d V_{1})$. Given $\Omega(d V_{1})=\Omega(d V_{2})$, we have
$\Omega(d U^{\bot}) \simeq \Omega(d U)$ on a removable set.
Note that given $p \leq \rho \leq q$, we have $\{ q < \lambda \} \subset \{ \rho < \lambda \} \subset \{ p < \lambda \}$,
that is we can assume $\rho \sim p_{1} \otimes p_{2}$. Note that a necessary condition for the inclusion
is that for instance $d U^{\bot} / d U = \rho \rightarrow 0$ in $\infty$, why $\rho(u,u^{\bot})$ 
essentially only depends on u. 

\newtheorem{orient}[prel]{Lemma}
\begin{orient}
Define $d N$ according to $d U \bot \sigma d I=d N$, with $p \leq \sigma \leq q$ 
for polynomials p,q, further $<q d I,d U>=<p d I,d U>=0$, we then have the orientation for $d U$ can be given relative a 
1-dimensional hyper plane.
\end{orient}
If instead $d N^{\bot}=d (I-N) + d V$, with $N^{\bot \bot} \sim N$, then $d V \sim d V^{\bot}$,
that the corresponding $\sigma$ are constants and we have a discontinuous continuation.
More generally, when $\frac{d N}{d N^{\bot}}=\frac{p}{q}$, where $p/q$ polynomial in $\infty$,
we have $d N^{\bot} \prec d N$ in $\infty$. When $d U= \rho d N$, a sufficient condition for 
independence of $d N^{\bot}$, is that $q \prec \prec p$ in $\infty$.
 \section{Riemann surfaces}
 \subsection{Planar Riemann surfaces}
 Riemann surfaces: define $\mathcal{W}$ through $d U=\rho d I$ on $\mathcal{W}$, implies $\rho=const$. We assume for Riemann surfaces
 that $\rho \in H(nbhd \ \mathcal{W})$ (analytic functions), that is removable sets are defined using continuation
 through continuity. When $\mathcal{W}$ Riemann, we assume a separation axiom of Hausdorff type.

Assume that $(u_{1},u_{2})$ defines a planar domain in $ L^{1}$, in particular
$L^{1}_{ac}$ with $\Omega(d U)=\Omega(U)$. Assume that $d V=\rho d U_{1}$, with 
$\rho \in L^{1}(u_{1},u_{2})$. We consider a continuation $L^{1}_{ac} \rightarrow L^{1}_{c}$.
Note that $L^{1}_{c}$ is equilibrated, as opposed to $L_{ac}^{1}$.
Example: given $ \sigma^{p} d L_{ac} = d L_{c}$, then $\sigma^{p} f^{p}=g \in L^{1}_{c}$. 
implies $\sigma f \in L_{c}^{1}$, but does not imply $\sigma f \in L_{ac}^{1}$.

The Dirichlet integral: given $u \in C^{1}(\mathcal{W})$ real, $D_{\mathcal{W}}(u)=\int \int_{\mathcal{W}} ( (\frac{\delta u}{\delta x})^{2} + (\frac{\delta u}{\delta y})^{2}) d x d y$.
We define HD: $\Delta u=0$ and $D_{\mathcal{W}}(u) < \infty$, we define AD: $u \in H(\mathcal{W})$ and $D_{\mathcal{W}}(u) < \infty$.
We write $\mathcal{W} \in O_{HD}$.
when $ u \in HD(\mathcal{W})$ implies $u=const$ and we write $\mathcal{W} \in O_{AD}$ when 
$ u$ $\in AD(\mathcal{W})$ implies $u=const$.
Obviously, $O_{HD} \subset O_{AD}$. (\cite{AhlforsSario60})
A parabolic Riemann surface $\mathcal{G}$, is defined through $\rho$ locally 1-1 and sub harmonic, implies 
$\rho=const$ on $\mathcal{G}$. Example: $Z \subset Z'$, zero spaces to $d U,d U'$ respective, where U is 
absolute continuous over $Z$ and $Z'$ parabolic, that is given
$d U' =\rho d U$ with $\rho$ locally 1-1 and sub harmonic over $Z'$, we then have
over $Z'$, that $U' f=f$.

Consider symmetric neighborhoods $(U_{1}-I, U_{2}) \rightarrow (U_{1}-I,-U_{2})$ and
$(I-U_{1},U_{2}) \rightarrow (U_{2},I-U_{1})$. Assume that $(I-U_{1})=U_{2} + V$, that is continuation
to $\mathcal{L}_{c}$ defines the polar. Assume that $\Phi_{ij} : L_{i} \rightarrow L_{j}$ a mapping 
between constant surfaces,
with $d U_{j}=\rho_{ij} d U_{i}$. Given $\overline{d} \rho_{ij}=0$
iff $\rho_{ij}=1$, on constant surfaces, defines $\Phi_{ij}$, as conjugation of invariant sets for 
$d U_{i},d U_{j}$.

Given leaves $\mathcal{L}_{j}$ and a transversal approximation T, we have T $\bot \mathcal{L}_{j}$
with unique orientation according to a normal model. Starting from a specific point $P \in \mathcal{L}_{j}$,
we require a movement on the leaves, to reach the transversal. Consider a transitive sub-group $H \subset G$, 
that preserves leaves, then the transversal can be reached from every P.   A spiral approximation involves two 
movements in the approximation, however we do not have a representation of bounded variation 
and not a unique orientation.

\subsection{Removable sets}
Given $\mathcal{W}$ planar and d U analytic
and bounded on $\mathcal{W}$, we then have $m (\complement \mathcal{W})=0$. Trivially,
given $\mathcal{W} \in \mathcal{O}_{AB}$ and $m(\complement \mathcal{W}) \neq 0$, we have $\mathcal{W}$ not planar, for instance $\rho,\vartheta \neq const$ i $\infty$,
with $d (U,U^{\bot})=(\rho d U_{1},\vartheta d U_{2})$.
Example: Consider the boundary 
$\{ \mid u_{3} \mid=0 \}$ and $\{ f \}$ $L^{1}-$ compact, in the sense that $m_{L^{1}} (\{ f(u) \quad \mid u_{3} \mid > R \}) \rightarrow 0$, 
when $R \rightarrow \infty$, that is $f_{u_{3}}$ planar $\rightarrow 0$ in $\infty$.

Assume that $d U(f)=-\{ h,f \}$, then we have existence of g so that $d U^{\diamondsuit}(f)=\{ g,f \}$,
with the involution condition $\xi^{\diamondsuit} / \eta^{\diamondsuit} = - \eta / \xi$.
Note that when $D_{\mathcal{W}}(h)=\int \xi^{2} + \eta^{2} d x d y$, where $\xi,\eta$ are polynomials, 
we have that $D_{\mathcal{W}}(h)=0$ implies $m(\mathcal{W})=0$.
Further, 
E is removable $\forall u \in HD$ iff $\complement E \in O_{G}$.
When $\mathcal{W}$ planar, we have $O_{G}=O_{HD}=O_{HB}=O_{HD}$

Example: consider $\mathcal{W}$ continued (simply connected) to $\infty$. 
For $\mathcal{O}_{HD}$ we have that $\mathcal{W}$ are first surfaces to harmonic functions, 
that is close to $\rho=const$, there is an orthogonal decomposition. Assume $\mathcal{W} \in \mathcal{O}_{AD}$, given $\alpha d U = d U_{1}$, 
$\beta d V= d U_{2}$ over $\tilde{\mathcal{W}} \ni \infty$ and $\overline{\delta} \alpha = \overline{\delta} \beta=0$, implies $\alpha,\beta = const$ on $\mathcal{W}$, 
we then have $\alpha=\beta=0$ on $\tilde{\mathcal{W}}$.
Over $\mathcal{W} \in \mathcal{O}_{G}$, given $d (U-I)=\rho d U$, that is $=\rho(u) d I$, where $\rho$ locally 1-1, we have
$d (U-I)f=0$ implies $f=0$.

\subsection{Normal movements}

Consider informally $\int f(u) d \varphi(u)$, where $f(u) \sim Uf \in (I)$, $d \varphi$ of bounded variation 
and $f(u)$ continuous and U analytic over f. 
Assume that $Uf \in (I)$ a compact space. The integral is indetermined over
$u=const$.  Assume that $G$ connected.
Further, assume that $d \varphi(u)$ monotonous, with $\int d \varphi(u)=1$, then the first mean value
theorem gives $\int f(u) d \varphi(u)=f(v)\int d \varphi \sim Vf$, that is we have existence of intermediary movements
$V \in G$, for instance $d V$ monotonous, with $U-I \subset V \subset U$.

Concerning normal operators, $\alpha$ a boundary for $\mathcal{W}' \subset \mathcal{W}$, regularly 
embedded in $\overline{\mathcal{W}}$ compact.
Assume that $\mathcal{W}'$ has a compact complement and boundary $\alpha$, in particular given $f \in L_{ac}^{1}$, 
we can assume $\alpha$ is defined by $U f(w)=f(u)=f(w)$ on $\alpha$. Assume that $f=f(n,u)$, where u=u(n), 
with $\frac{d f}{d n}$ bounded. Then $\beta$ can be defined by $\frac{\delta f}{\delta n}=0$ on $\beta$. 
Note that $f=f(u)$ on $\beta$ and we can assume for instance $N(f)=1$ on $\beta$. According to the theory 
on evolute, if $\Sigma=\{ d U=0 \}$ and $T \Sigma=\{ \Delta U=0 \}$,
we have $(T \Sigma)^{\bot}=T \Sigma^{\bot}$, where we can define $\Sigma^{\bot}(\varphi)=\{\varphi \quad d N \bot d U \}$,
that is $<d N, d U>(\varphi)=0=UN(\varphi)$, where $N(\varphi) \neq 0$. 
Where $d U + d N=d I$, over $\Omega(d U)$, we have $N(\phi)=0$.
When $N=N(U_{1},U_{2})$, since $d U_{1} \bot d U_{2}$ over a (linearly) convex set, we must have $d N / d U_{2} \prec \prec d N/d U_{1}$.
Continue $f$ to $L_{c}^{1}$, that is compact disk neighborhoods. Given $N$ harmonic, we have that 
$N \rightarrow {}^{t} N$ preserves type, that is $(N,{}^{t} N)$ is symmetric. There is analogously 
with Ahlfors, a normal id-operator
for $R(G) \cap L_{c}^{1}$, under the condition that $R(G)^{\bot}$ $L^{1}-$ compact, with an approximation property.

\section{Spirals}

\subsection{Non-regular approximations}
Assume that $N(e^{\phi})=e^{n(\phi)}$, given n single valued over a domain of $\phi$,
assume that $n(\phi)=<g,\phi>$. Note that $\phi \in \mathcal{D}_{L^{1}}$ can be approximated by
polynomials and assume n linear over $\phi$, where $\phi \in R(U_{1},U_{2})$ is algebraic in 
$u_{1}$ $(u_{2})$. Given $d n=\rho_{1}d u_{1}$, with $\rho_{j} \in AD$ single valued, we have that
N algebraic in $u_{1}$ $(u_{2})$ over the domain. Over a strictly pseudo convex domain, we have that
N has locally algebraic definition and there is not space for a spiral.
Assume that $n(\phi) \in L^{1}$, that is $(u_{1},u_{2})$ dense in the domain, $n(\phi)=g*\phi(u_{1})=\int \phi(u_{1} - u_{2}) d n(u_{2})$,
continue N to $ \tilde{N}$: $n(\phi) + v(u_{1},u_{2})$, where locally $v(u_{1},u_{2})=1$ over a domain and $v \in \mathcal{H}$ over the domain. The condition 
$\tilde{N}(e^{\phi})=e^{\phi}$ under the condition $g \bot \phi$
means $v - I(\phi)=0$ iff $1 = \phi(u_{1}-u_{2})$, where we assume $\phi(0)=1$, that is vi consider
$\phi \in \mathcal{D}_{L^{1}}$, with continuation through constant surfaces. In this case
there is space for a spiral in the topology $\mathcal{L}_{c}$, through the diagonal $u_{1}=u_{2}$.
More precisely, $f(x,y) \rightarrow f(u_{1},u_{2})$, since $\frac{\delta f}{\delta u_{1}} \frac{\delta u_{1}}{\delta x} + \frac{\delta f}{\delta u_{2}} \frac{\delta u_{2}}{\delta x}$,
 given $f(u_{1},u_{2}) \equiv 1$ we have, when f absolute continuous in x,y, $f(x,y) \equiv const$. However, $T(u,v) - \delta_{0}(v-u) \sim 0$,
 does not necessarily imply $T(x,y) - \delta_{x}(y) \sim 0$.

 Assume the Lie group corresponding to
$E'$ has maximal order, that is that the domain that is the result of change of variables, is dense in 
the domain, for $T \in \mathcal{H}$.
Assume that $F : L_{ac}^{1} \rightarrow L_{ac}^{1}$ with $F(x,y)-Id(x,y)=Q(x,y)$, where $Q$ polynomial 
(approximates $C^{\infty}$).
When $Q(x,y) \leq F(x,y)$, then F has semi-algebraic sub level surfaces. Since Id can be given 
injectively on $L_{ac}^{1}$ and $Q$ has removable zero's, we can consider F as projective on $L_{ac}^{1}$. 
For $F : L_{c}^{1} \rightarrow L_{c}^{1}$, Id has support on the diagonal, which leaves space for 
non-trivial kernels for F. The same applies, when we continue
Q to F very regular i $L_{c}^{1}$.

Note that starting from en planar domain $(U,U^{\bot})$ and $\Gamma =\{ U=U^{\bot} \}$, we can map 
a closed Jordan arc on $\Gamma$, through the logarithm mapping. Given $d U,d U^{\bot}$ has constant 
opposite sign, there is corresponding to the planar domain, a corresponding planar Riemann 
domain. 
  
 Assume that $d U(f)=-\{ u,f \}$ and $d V(f)=-\{ v,f \}$, with u,v analytic. Given u,v conjugated 
 according to
 $-i(\frac{\delta u}{\delta y},-\frac{\delta u}{\delta x}) = (-\frac{\delta v}{\delta y},\frac{\delta v}{\delta x})$,
 analyticity is preserved for the corresponding forms. Consider continuation along $L=\{ (x,y) \}$
 to $\infty$. Assume that $d V(x^{*},y^{*}) =\sigma(x,y) d U(x,y)$, where $x x^{*} + y y^{*}=R \neq 0$,
 that is a contact transform. Given $d V_{j}$ with linear independence i $\infty$, then where $\sigma=const$, 
 the corresponding $d U_{j}$ are linearly
 independent close to 0. We assume $\mid d U(x,y) - d V(\frac{1}{x},\frac{1}{y}) \mid \leq \epsilon$ 
 for $(x,y)$ close to $\infty$. Note that given the
 approximation property for R(d V), we have that d U according to the above, preserves constant value 
 in the $\infty$.
 
 Example: assume that $U+V=I$ projective in phase, we then have given $R(UV)$ has the approximation
 property, over analytic f, that $UVf(x,y) \rightarrow f(u,1/u)$ is continuous. Consider $<E(x,y),\phi \otimes \psi>$ 
 and $V^{-1}U \phi=\phi$,$U^{-1}V \psi=\psi$, we then have 
 $<E(u,v),\phi \otimes \psi>=<E(v,u),\phi \otimes \psi>$
 
Assume that $d U=\rho d I$,
with $\rho \in AD$.
Given $\rho$ non linear, over a planar domain, we have that $d U$ is not single valued. Single valuedness concerns the coefficients 
to d U. 
Consider $(d U)^{\bot}$ multivalent, given holomorphic leaves, the spiral gives a (non-regular) approximation 
property. Consider $\Psi : d U \rightarrow d U^{\bot}$, where $\Psi$ sub nuclear in the tangent space, 
note given that not both movements are absolute continuous, we do not necessarily have $\mid U^{\bot} - I \mid \leq \mid U-I \mid$. 
However, $\mid U^{\bot} -  U \mid \leq \mid U^{\bot} - I \mid + \mid U-I \mid$ and
over functions symmetric with respect to the right hand sides, for which the movements are analytic, 
we can present a spiral.

\bibliographystyle{amsplain}
\bibliography{sos}
\end{document}